\documentclass[amsmath,10pt,a4paper]{article}
\usepackage{amsfonts, amssymb, amsmath, amscd}
\usepackage[french]{babel}
\usepackage{amssymb}
\newtheorem{thm}{Th\'eor\`eme}[section]
\newtheorem{conj}[thm]{Conjecture}
\newtheorem{lem}[thm]{Lemme}
\newtheorem{prop}[thm]{Proposition}
\newtheorem{cor}[thm]{Corollaire}

\newtheorem{Q}[thm]{Question}
\newtheorem{prob}[thm]{Probl\`eme}

\numberwithin{equation}{section}

\title{Repr\'esentations cohomologiques isol\'ees, applications cohomologiques}

\author{N. Bergeron}

\date{}

\begin{document}

\maketitle

\section*{Introduction}

Un groupe de Lie $G$ a la propri\'et\'e (T) de Kazhdan si la repr\'esentation triviale de $G$ est isol\'ee dans le dual unitaire 
$\widehat G$ de $G$, muni de sa topologie de Fell. 
D'apr\`es Kazhdan \cite{Kazhdan} et Kostant \cite{Kostant}, un groupe semi-simple r\'eel $G$ a la propri\'et\'e (T)
si et seulement s'il est sans facteur simple localement isomorphe \`a $SO(n,1)$ ($n\geq 2$) ou $SU(n,1)$ ($n\geq 1$). La repr\'esentation triviale du groupe 
$G=SO(n,1)$, $n\geq 2$ (resp. $G=SU(n,1)$, $n\geq 1$) n'est en particulier pas isol\'ee dans le dual unitaire de $G$. Selberg montre n\'eanmoins dans \cite{Selberg} qu'il subsiste une 
propri\'et\'e d'isolation, de nature arithm\'etique, dans le cas du groupe $G=SL(2,{\Bbb R})$ (localement isomorphe au groupe $SO(2,1)$)~: la repr\'esentation triviale 
du groupe $G$ est isol\'ee dans le sous-ensemble de $\widehat G$ constitu\'e de toutes les repr\'esentations apparaissant faiblement dans l'une des repr\'esentations
quasi-r\'eguli\`eres $L^2 (\Gamma (N) \backslash G)$ avec $N \in {\Bbb N}$ et $\Gamma (N) = {\rm ker} (SL(2,{\Bbb Z}) \rightarrow SL(2,{\Bbb Z} /N{\Bbb Z} ))$.

Fixons plus g\'en\'eralement $G$ une groupe semi-simple alg\'ebrique sur ${\Bbb Q}$. Soit $G({\Bbb R})^0$ la composante neutre de $G({\Bbb R})$~: c'est un produit de groupes 
semi-simples r\'eels. Nous notons 
$$G({\Bbb R})^0 = U \times G^{{\rm nc}},$$
o\`u $U$ et compact et $G^{{\rm nc}}$ est un groupe semi-simple r\'eel sans facteur compact. Un sous-groupe $\Gamma \subset G({\Bbb Q})$ est dit {\it de congruence} s'il existe 
un sous-groupe compact-ouvert $K_f$ de $G({\Bbb A}_f )$ (o\`u ${\Bbb A}_f$ d\'esigne l'anneau des ad\`eles finis sur ${\Bbb Q}$) tel que $\Gamma = G({\Bbb Q}) \cap K_f$. Suivant 
Burger et Sarnak \cite{BurgerSarnak}, on appelle {\it spectre} de $\Gamma \backslash G$ - o\`u $\Gamma$ est un sous-groupe de congruence de $G$ - l'ensemble des repr\'esentations
irr\'eductibles unitaires $\pi \in \widehat{G}$ apparaissant (faiblement, si $G$ est isotrope sur ${\Bbb Q}$) dans la repr\'esentation r\'eguli\`ere de $G$ dans $L^2 (\Gamma \backslash G)$~:
\begin{eqnarray} \label{spectre}
\sigma (\Gamma \backslash G) = \{ \pi \in \widehat{G} \; : \; \pi \propto L^2 (\Gamma \backslash G) \} . 
\end{eqnarray}
Rappelons \'egalement la d\'efinition du dual automorphe $\widehat{G}_{{\rm Aut}}$ de $G$ donn\'ee dans \cite{BurgerSarnak}~:
\begin{eqnarray} \label{dualautom}
\widehat{G}_{{\rm Aut}} = \overline{\cup_{\Gamma \; {\rm cong}} \sigma (\Gamma \backslash G) }.
\end{eqnarray}

Le Th\'eor\`eme de Selberg mentionn\'e plus haut affirme alors que si $G= SL(2)_{|{\Bbb Q}}$, la repr\'esentation trivial $1_G$ de $G$ est {\bf isol\'ee} dans le dual automorphe 
$\widehat{G}_{{\rm Aut}}$. Cette propri\'et\'e, appel\'e {\it propri\'et\'e $\tau$} est en fait vraie pour tout groupe $G$ comme ci-dessus, cf. la d\'emonstration de la Conjecture $\tau$ par Clozel \cite{Clozel}.

Dans cet article nous \'etudions des g\'en\'eralisations de la propri\'et\'e (T) et de la propri\'et\'e $\tau$. La repr\'esentation triviale de $G$ est l'exemple le plus simple d'une {\bf repr\'esentation 
cohomologique}, c'est-\`a-dire (essentiellement, cf. \S 1 pour une d\'efinition pr\'ecise) d'une repr\'esentation pouvant intervenir dans la cohomologie d'une vari\'et\'e localement sym\'etrique. La repr\'esentation
triviale est cohomologique de degr\'e $0$, elle intervient dans les fonctions constantes. Les g\'en\'eralisations des propri\'et\'es (T) et $\tau$ que nous avons en vus ont traits aux deux probl\`emes suivants~:

\begin{prob} \label{pbm1}
D\'ecrire les repr\'esentations isol\'ees dans $\widehat G$ parmi les 
repr\'esenta\-tions cohomologiques.
\end{prob}

\begin{prob} \label{pbm2}
Pour $G/{\Bbb Q}$ donn\'e, soit $\pi \in \widehat G$ une 
repr\'esentation cohomologique non triviale et non isol\'ee. La
repr\'esentation $\pi$ est-elle isol\'ee {\bf dans} $\widehat G_{\rm Aut} \cup \{ \pi \}$ ?
\end{prob}

Remarquons imm\'ediatement qu'il existe une classe simple de repr\'esentations pour lesquelles les deux questions ont une solution n\'egative~: la classe des repr\'esentations cohomologiques temp\'er\'ees, cf. \cite{Livre}.
Le but de cet article est d'apporter des r\'eponses tant\^ot conjecturale tant\^ot inconditionnelles \`a ces deux probl\`emes. Le Probl\`eme \ref{pbm1} a en fait d\'ej\`a \'et\'e compl\`etement r\'esolu par Vogan, dans 
un article longtemps clandestin \cite{Vogan}. Sa solution est difficile, nous en donnons, au \S 1 une d\'emonstration ``\'el\'ementaire'' au sens o\`u nous n'utilisons ``que'' la classification, par Vogan et Zuckerman 
\cite{VoganZuckerman}, \cite{Vogan2}, des repr\'esentations cohomologiques.

Le Probl\`eme \ref{pbm2} est encore plus profond~: l'un des buts de cet article (cf. \S 3) est d'expliquer comment sa solution est implicitement donn\'ee par les Conjectures d'Arthur. Nous montrons au \S 4 que ce
Probl\`eme n'est en tout cas pas gratuit en obtenant un certain nombre de cons\'equences cohomologiques, pour certaines nouvelles, concernant les vari\'et\'es arithm\'etiques. 

Illustrons imm\'ediatement nos r\'esultats par un exemple. Soit $G$ un groupe semi-simple alg\'ebrique sur ${\Bbb Q}$ tel que $G^{{\rm nc}} = O(p,q)^0$. Supposons de plus que 
$G$ ne provient pas d'un groupe de type ${}^{3,6} D_4$. Soit $\mathfrak{q}$ une sous-alg\`ebre parabolique de l'alg\`ebre de Lie complexifi\'ee de $O(p,q)$, stable par l'involution de Cartan. Soit $L$ la composante 
connexe du normalisateur de $\mathfrak{q}$ dans $O(p,q)^0$. Supposons qu'il existe un entier naturel $r$, $2r \leq q$, tel que
\begin{eqnarray} \label{Lintro}
L \cong U(1)^{r} \times O (p, q-2r)^0 .
\end{eqnarray}
Dans \cite{VoganZuckerman} Vogan et Zuckerman associent \`a $\mathfrak{q}$ une repr\'esentation cohomologique $A_{\mathfrak{q}}$ intervenant en degr\'e $R=pr$. (Les r\^oles
de $p$ et $q$ pourraient naturellement \^etre \'echang\'es.)

Le th\'eor\`eme de Vogan (r\'eponse au Probl\`eme \ref{pbm1}) dont nous donnons une d\'emonstration ``\'el\'ementaire'' au \S 1 implique le th\'eor\`eme suivant.

\begin{thm} \label{T1intro}
La repr\'esentation cohomologique $A_{\mathfrak{q}}$ est {\bf isol\'ee} dans le dual unitaire du groupe $O(p,q)^0$ si et seulement si
$p\geq 2$, $q\geq 2r+2$ et $p+q \geq 2r+5$.
\end{thm}

Concernant le Probl\`eme \ref{pbm2}, la Conjecture \ref{isolaut} que nous motivons au \S 3  implique la conjecture suivante.

\begin{conj} \label{C1intro}
Pour tout entier naturel $r \leq q/2$, la repr\'esentation cohomologique $A_{\mathfrak{q}}$ est {\bf isol\'ee} dans $\widehat G_{\rm Aut} \cup \{A_{ \mathfrak{q}} \}$.
\end{conj}

Ces consid\'erations nous permettent de montrer le th\'eor\`eme suivant au \S 4.

\begin{thm} \label{T2intro}
Supposons la repr\'esentation $A_{\mathfrak{q}}$ isol\'ee dans $\{ A_{\mathfrak{q}} \} \cup \widehat{G}_{{\rm Aut}}$
(c'est en particulier le cas si $p\geq 2$, $g\geq 2r+2$ et $p+q \geq 2r+5$, cf. Thm. \ref{T1intro}). Il existe alors un sous-groupe de congruence $\Gamma \subset G({\Bbb Q})$ tel que la repr\'esentation $A_{\mathfrak{q}}$ intervienne discr\`etement
dans $L^2 (\Gamma \backslash G)$.
\end{thm}

Notons $X_{p,q}$ l'espace sym\'etrique associ\'e au groupe $O(p,q)^0$. Supposons $p\leq q$. \`A l'aide de la formule de Matsushima (cf. \S 1) le Th\'eor\`eme ci-dessus (et un peu plus lorsque $p=1$) implique(nt)~:

\begin{thm} \label{T3intro}
Pour tout sous-groupe de congruence $\Gamma$ suffisamment profond, le $p$-i\`eme nombre de Betti de l'espace localement sym\'etrique
$\Gamma \backslash X_{p,q}$ est non nul.
\end{thm}

Il y a une litt\'erature abondante sur ce type de r\'esultats le Th\'eor\`eme \ref{T2intro} est en particulier d\^u \`a Li \cite{Li2} lorsque $r<\frac{1}{4} (p+q-2)$, il en d\'eduit le Th\'eor\`eme \ref{T3intro}
pour $p+q \geq 7$. Les premiers r\'esultats de ce type ont \'et\'e obtenus par Millson dans le cas du groupe $O(1,n)$.

\medskip

Pour finir, nous montrons, au \S 5, un r\'esultat, dont nous pensons qu'il est int\'eressant par lui-m\^eme, qui permet de ramener, dans un grand nombre de cas,
le Probl\`eme \ref{pbm2} au cas o\`u le groupe $G^{\rm nc} = O(n,1)^0$ ou $SU(n,1)$, tous deux explor\'es en d\'etails dans \cite{Livre}.

\section{Un th\'eor\`eme de Vogan}

Dans cette section $G=G^{\rm nc}$. Le groupe $G$ est donc semi-simple r\'eel et connexe. Fixons $K$ un sous-groupe compact maximal de $G$.
Notons $\mathfrak{g}_0$ l'alg\`ebre de Lie de $G$ et $\mathfrak{g}_0 = \mathfrak{k}_0 \oplus \mathfrak{p}_0$ la
d\'ecomposition associ\'ee au choix de $K$. Si $\mathfrak{l}_0$ est une alg\`ebre de Lie, nous noterons
$\mathfrak{l}=\mathfrak{l}_0 \otimes {\Bbb C}$ sa complexification.

Rappelons qu'un {\it $(\mathfrak{g} , K)$-module} est un espace vectoriel complexe $V$ muni de repr\'esentations de 
$\mathfrak{g}$ et $K$, soumises aux trois conditions suivantes.
\begin{enumerate}
\item L'action de $K$ sur $V$ est {\it localement finie}, {\it i.e.} tout vecteur $v \in V$ appartient \`a un sous-espace
$K$-invariant $V_1 \subset V$ de dimension finie, et la repr\'esentation de $K$ dans $V_1$ est lisse.
\item La diff\'erentielle de l'action de $K$ (qui est bien d\'efinie d'apr\`es 1.) est \'egale \`a la restriction \`a $\mathfrak{k}$
de l'action de $\mathfrak{g}$.
\item Si $X \in \mathfrak{g}$, $k\in K$ et $v \in V$, alors $k \cdot (X \cdot v) = ({\rm Ad} (k) X) \cdot (k \cdot v)$.
\end{enumerate} 
Si $(\pi , {\cal H}_{\pi} )$ est un repr\'esentation continue de $G$ dans un espace de Hilbert, l'espace des vecteurs lisses et 
$K$-finis de $\pi$ est un $(\mathfrak{g} , K)$-module appel\'e {\it module de Harish-Chandra de $\pi$}.

Soit maintenant $(\pi , V_{\pi} )$ un $(\mathfrak{g} , K)$-module irr\'eductible. \`A tout r\'eseau $\Gamma \subset G$, on peut naturellement associer
une application lin\'eaire
\begin{eqnarray} \label{Tpi}
T_{\pi} (\Gamma ) : \left\{
\begin{array}{rcl}
\mbox{Hom}_K (\bigwedge^* \mathfrak{p}, \pi ) \otimes \mbox{Hom}_{\mathfrak{g} , K} (\pi , C^{\infty} (\Gamma \backslash G ))  & \rightarrow & {\cal E}^* (S(\Gamma )) , \\
\psi \otimes \varphi & \mapsto & \varphi \circ \psi ,
\end{array} \right.
\end{eqnarray}
\`a valeurs dans l'espace ${\cal E}^* (S(\Gamma ))$ des formes diff\'erentielles sur $S(\Gamma ) = \Gamma \backslash G /K$.
L'application $T_{\pi}(\Gamma )$ associe donc \`a certains $(\mathfrak{g}, K)$-modules $\subset C^{\infty} (\Gamma \backslash G)$ des formes diff\'erentielles sur $S(\Gamma )$, l'action
du laplacien de Hodge-de Rham sur les formes diff\'erentielles est induit par l'action de  l'{\it op\'erateur de Casimir} 
$$\Omega = \sum_{1\leq s \leq n} y_s . y_s '$$
o\`u $(y_s)$ est une base de $\mathfrak{g}$ et $(y_s ')$ la base duale de $\mathfrak{g}$ par rapport
\` a sa forme de Killing, sur $\pi$. Un $(\mathfrak{g} , K)$-module ne peut donc \'eventuellement intervenir dans la cohomologie d'une vari\'et\'e localement sym\'etrique 
$S(\Gamma )$ que s'il est unitaire et v\'erifie
\begin{enumerate}
\item $\pi (\Omega ) = 0$, et
\item $\mbox{Hom}_K (\bigwedge^* \mathfrak{p}, \pi ) \neq \{ 0 \}$.
\end{enumerate}
Un $(\mathfrak{g} , K)$-module v\'erifiant les points 1. et 2. ci-dessus est dit {\it cohomologique}, il est \'equivalent de dire que sa $(\mathfrak{g} , K)$-cohomologie est non triviale.
Rappelons que la {\it $(\mathfrak{g} , K)$-cohomologie de $\pi$}~: $H^* ( \mathfrak{g} , K ; \pi )$, est la cohomologie du complexe 
$$C^* (\mathfrak{g} , K ; \pi ) := {\rm Hom}_K (\bigwedge {}^* \mathfrak{p} , \pi ),$$
muni de la diff\'erentielle d\'efinie pour $\omega \in C^p  (\mathfrak{g} , K ; \pi )$ par la formule~:
$$d\omega (X_0 , \ldots , X_p )  = \sum_{i=0}^p (-1)^i X_i \cdot \omega (X_0 , \ldots , \widehat{X_i} , \ldots , X_p ) .$$

D'apr\`es Parthasarathy \cite{Parthasarathy}, Kumaresan \cite{Kumaresan} et Vogan-Zuckerman \cite{VoganZuckerman},
les repr\'esentations cohomologiques sont associ\'ees aux sous-alg\`ebres paraboliques $\theta$-stable de $\mathfrak{g}$. Rappelons en la d\'efinition.
Notons $\mathfrak{t}_0 = $Lie$(T)$ une sous-alg\`ebre de Cartan de $\mathfrak{k}_0$. 
Une {\it sous-alg\`ebre parabolique $\theta$-stable} $\mathfrak{q} = \mathfrak{q} (X) \subset \mathfrak{g}$ est associ\'ee \`a un 
\'el\'ement $X\in i \mathfrak{t}_0$. Elle est \'egale \`a la somme 
$$\mathfrak{q}  = \mathfrak{l}  \oplus \mathfrak{u},$$ 
du centralisateur $\mathfrak{l}$ de $X$ et du sous-espace $\mathfrak{u}$ engendr\'e par les racines positives de $X$ dans $\mathfrak{g}$.
Alors $\mathfrak{q} $ est stable sous $\theta$; on en d\'eduit une d\'ecomposition $\mathfrak{u}  = (\mathfrak{u}  \cap
\mathfrak{k} ) \oplus (\mathfrak{u}  \cap \mathfrak{p} )$. Soit $R =
\mbox{dim} (\mathfrak{u} \cap \mathfrak{p} )$ ({\it degr\'e fortement primitif}).

Associ\'e \`a $\mathfrak{q}$, se trouve un $(\mathfrak{g} , K)$-module irr\'eductible bien d\'efini $A_{\mathfrak{q}}$ que nous caract\'erisons
maintenant. Supposons effectu\'e un choix de racines positives pour $(\mathfrak{k} , \mathfrak{t} )$ de fa\c{c}on
compatible avec $\mathfrak{u}$. Soit $e(\mathfrak{q})$ un g\'en\'erateur de la droite $\bigwedge^R (\mathfrak{u} \cap \mathfrak{p} )$.
Alors $e(\mathfrak{q})$ est le vecteur de plus haut poids d'une repr\'esentation irr\'eductible $V(\mathfrak{q})$ de $K$ contenue
dans $\bigwedge^R \mathfrak{p}$; et dont le plus haut poids est donc n\'ecessairement $2\rho (\mathfrak{u} \cap \mathfrak{p} )$.
La classe d'\'equivalence du $(\mathfrak{g} , K)$-module $A_{\mathfrak{q}}$ est alors uniquement caract\'eris\'ee
par les deux propri\'et\'es suivantes.
\begin{eqnarray} \label{Aq1}
\begin{array}{l}
A_{\mathfrak{q}} \mbox{ {\it est unitarisable avec le m\^eme caract\`ere infinit\'esimal que la}} \\
\mbox{{\it  repr\'esentation triviale}}
\end{array}
\end{eqnarray}
\begin{eqnarray} \label{Aq2}
\mbox{Hom}_K (V(\mathfrak{q} ), A_{\mathfrak{q}} ) \neq 0.
\end{eqnarray}
Remarquons que la classe du module $A_{\mathfrak{q}}$ ne d\'epend que de l'intersection $\mathfrak{u} \cap \mathfrak{p}$,
autrement dit deux sous-alg\`ebres paraboliques $\mathfrak{q} = \mathfrak{l} \oplus \mathfrak{u}$ et $\mathfrak{q}' =\mathfrak{l}' \oplus
\mathfrak{u}'$ v\'erifiant $\mathfrak{u} \cap \mathfrak{p} = \mathfrak{u}' \cap \mathfrak{p}$ donnent lieu \`a une m\^eme classe de
module cohomologique.

De plus, $V(\mathfrak{q})$ intervient avec multiplici\'e $1$ dans $A_{\mathfrak{q}}$ et $\bigwedge^R (\mathfrak{p} )$, et
\begin{eqnarray} \label{gKcohom}
H^i (\mathfrak{g} , K , A_{\mathfrak{q}} ) \cong \mbox{Hom}_{L\cap K} ( \bigwedge {}^{i-R} (\mathfrak{l} \cap \mathfrak{p} ), {\Bbb C}).
\end{eqnarray}
Ici $L$ est un sous-groupe de $K$ d'alg\`ebre de Lie $\mathfrak{l}$.

Le th\'eor\`eme fondamental suivant, d\^u \`a Vogan et Zuckerman d\'ecoule alors de \cite{VoganZuckerman}.

\begin{thm} \label{rep cohom}
Soit $\pi$ un $(\mathfrak{g} , K)$-module unitaire cohomologique. Il existe une sous-alg\`ebre parabolique
$\theta$-stable $\mathfrak{q}=\mathfrak{l}+\mathfrak{u}$ de $\mathfrak{g}$ tel que
$\pi \cong A_{\mathfrak{q}}$. On peut de plus choisir $\mathfrak{q}$ de mani\`ere \`a ce que
le groupe $L$ n'ait pas de facteurs simples compacts (non ab\'elien); dans ce cas
$\mathfrak{q}$ est uniquement d\'etermin\'e \`a conjugaison par $K$ pr\`es.
\end{thm}

Nous supposerons dor\'enavant que l'alg\`ebre $\mathfrak{q}$ associ\'ee \`a une repr\'esentation
$\pi$ comme dans le Th\'eor\`eme \ref{rep cohom} est choisie de mani\`ere \`a ce que
le groupe $L$ n'ait pas de facteurs simples compacts (non ab\'elien).

Nous allons maintenant d\'eduire du Th\'eor\`eme \ref{rep cohom}, le th\'eor\`eme suivant
qui est un cas particulier d'un th\'eor\`eme \cite[Theorem A.10]{Vogan} de Vogan. Notons $\widehat{G}$ le {\it dual unitaire} de $G$, d'apr\`es un th\'eor\`eme classique d'Harish-Chandra, celui-ci
s'identifie \`a l'ensemble des classes d'\'equivalences
de $(\mathfrak{g}, K)$-modules unitaires irr\'eductibles, nous le munissons de sa topologie de Fell. 
Notons enfin $\Pi$ (resp. $\Pi (\mathfrak{l})$) l'ensemble des racines simples de $\mathfrak{t}$ dans $\mathfrak{g}$ (resp. $\mathfrak{l}$).

\begin{thm} \label{rep isol}
Soit $\pi$ comme dans le Th\'eor\`eme \ref{rep cohom}. La repr\'esentation $\pi$ est {\bf isol\'ee}
dans le dual unitaire de $G$ si et seulement si le couple $(G , \mathfrak{q})$ a les propri\'et\'es
suivantes.
\begin{enumerate}
\item Le groupe $L$ n'a aucun facteur simple localement isomorphe \`a $SO(n,1)$ ou $SU(n,1)$ ($n\geq 1$).
\item Il n'y a pas de racine imaginaire non compacte dans $\Pi$ orthogonale \`a $\Pi (\mathfrak{l})$.
\end{enumerate}
\end{thm}
{\it D\'emonstration.} Nous commen\c{c}ons par d\'emontrer que les conditions 1. et 2. sont suffisantes. Nous montrons en fait la contrapos\'ee. Supposons donc que la
repr\'esentation $\pi$ n'est pas isol\'ee dans le dual unitaire de $G$~: soit $\{ \pi_j \}$ une suite de repr\'esentations irr\'eductibles unitaires diff\'erentes de
$\pi$ qui converge vers $\pi$. 

Le plus bas $K$-type de $\pi$ est un $K$-type des repr\'esentations $\pi_j$ (sauf
peut-\^etre un nombre fini d'entre elles). Notons alors $\mu = \mu (\mathfrak{q})$
le plus bas $K$-type de $\pi$, $R = R(\mathfrak{q} )= {\rm dim} (\mathfrak{u} \cap \mathfrak{p} )$ et supposons dor\'enavant
que $\mu$ intervient comme $K$-type dans chaque $\pi_j$. Le $K$-type $\mu$ intervient avec multiplicit\'e $1$ dans $\pi$
et $\bigwedge^R \mathfrak{p}$, et
$$H^* (\mathfrak{g} , K , \pi ) \cong {\rm Hom}_K (\bigwedge {}^* \mathfrak{p},  \pi ) \cong {\rm Hom}_{L\cap K} (\bigwedge {}^{*-R} (\mathfrak{l} \cap \mathfrak{p}) , {\Bbb C}) .$$

Fixons un entier $j$ et consid\'erons le complexe
$$\ldots \rightarrow C^* (\pi_j ) := {\rm Hom}_K (\bigwedge {}^* \mathfrak{p} , \pi_j  ) \stackrel{d}{\rightarrow} C^{*+1} (\pi_j ) \rightarrow \ldots $$
Rappelons que si $f \in C^q (\pi_j )$, la diff\'erentielle $df$ de $f$ est l'\'el\'ement de $C^{q+1} (\pi_j )= {\rm Hom}_K (\bigwedge^{q+1} \mathfrak{p} , \pi_j  ) $ donn\'e par
$$df (x_0 , \ldots , x_q ) = \sum_i (-1)^i x_i \cdot f(x_0 , \ldots , \widehat{x_i} , \ldots , x_q ).$$

Puisque $\mu \subset \bigwedge^R \mathfrak{p}$, le produit scalaire sur $\bigwedge^* \mathfrak{p}$ induit par la forme de
Killing de $G$ permet d'identifier canoniquement ${\rm Hom}_K (\mu , \pi_j )$ \`a un sous-espace de
$C^R (\pi_j ) = {\rm Hom}_K (\bigwedge {}^R \mathfrak{p} , \pi_j  )$. De la m\^eme mani\`ere
l'inclusion $\mathfrak{p} \wedge \mu \subset \bigwedge^{R+1} \mathfrak{p}$ induit une inclusion
${\rm Hom}_K (\mathfrak{p} \wedge \mu , \pi_j ) \subset C^{R+1} (\pi_j )$. Et, si $f$ est un \'el\'ement non nul de
${\rm Hom}_K (\mu , \pi_j )$, par d\'efinition de la diff\'erentielle $d$,
l'\'el\'ement $df \in C^{R+1} (\pi_j )$ appartient en fait au sous-espace ${\rm Hom}_K (\mathfrak{p} \wedge \mu , \pi_j )$.

Nous distingons alors deux cas suivant que l'image $d({\rm Hom}_K (\mu , \pi_j ))$ soit triviale ou non \`a partir
d'un certain rang $j \geq j_0$.

\bigskip

Commen\c{c}ons par supposer que l'image $d({\rm Hom}_K (\mu , \pi_j ))$ est non triviale pour une infinit\'e de $j$.
Quitte \`a extraire on peut alors supposer $d({\rm Hom}_K (\mu , \pi_j )) \neq 0$ pour tout $j$. Comme repr\'esentation
de $K$, le module $\mathfrak{p} \wedge \mu$ se d\'ecompose en un nombre fini de sous-espaces irr\'eductibles. Quitte
\`a extraire une sous-suite de $(\pi_j )$, on peut donc supposer qu'il existe un $K$-type $\nu \subset \mathfrak{p}
\wedge \mu$ qui intervient dans chaque $\pi_j$. Pour conclure dans ce cas rappelons la c\'el\`ebre in\'egalit\'e de
Dirac de Parthasarathy (cf. \cite[(2.26)]{Parthasarathy}, \cite[II.6.11]{BorelWallach}, \cite[Lemma 4.2]{VoganZuckerman}).

\begin{lem} \label{dirac}
Soit $(\pi, V_{\pi} )$ une repr\'esentation unitaire irr\'eductible de $G$ et $V_{\pi}^K$ son $(\mathfrak{g} , K)$-module
associ\'e. Fixons une repr\'esentation de $\mathfrak{k}$ de plus haut poids $\chi \in \mathfrak{t}^*$ et apparaissant
dans $V_{\pi}^K$; et un sous-syst\`eme positif de racines $\Delta^+ (\mathfrak{g} )$ de $\mathfrak{t}$ dans
$\mathfrak{g}$. Notons $\rho$ (resp. $\rho_c$, $\rho_n$) dans $\mathfrak{t}^*$ la demi-somme des racines dans
$\Delta^+ (\mathfrak{g} )$ (resp. $\Delta^+ (\mathfrak{k} )$, $\Delta^+ (\mathfrak{p})$), de sorte que $\rho=\rho_c
+ \rho_n$. Soit $w$ un \'el\'ement du groupe de Weyl $W_K = W(\mathfrak{k} , \mathfrak{t} )$ de $\mathfrak{t}$ dans
$\mathfrak{k}$, tel que $w(\chi - \rho_n )$ soit dominant pour $\Delta^+ (\mathfrak{k})$. Alors,
$$-\pi (\Omega ) \geq ||\rho ||^2 -||w(\chi - \rho_n ) + \rho_c ||^2,$$
o\`u $\Omega$ d\'esigne le casimir de $\mathfrak{g}$ et la norme $||.||$ est d\'eduite de la forme de Killing sur
$\mathfrak{g}$.
\end{lem}

Puisque $\pi (\Omega )=0$ ($\pi$ est cohomologique et donc de caract\`ere infinit\'esimal \'egal \`a celui de la
repr\'esentation triviale), la suite $\pi_j (\Omega )$ tend vers $0$ lorsque $j$ tend vers $+\infty$. Le Lemme \ref{dirac}
appliqu\'e \`a chacune des repr\'esentations $\pi_j$ et au $K$-type $\nu$ implique que le plus haut poids $\chi$
de $\nu$ v\'erifie l'in\'egalit\'e
\begin{eqnarray} \label{eg}
||\rho ||^2 -||w(\chi - \rho_n ) + \rho_c ||^2 \leq 0.
\end{eqnarray}
Par ailleurs Kumaresan montre dans \cite{Kumaresan} qu'un $K$-type de $\bigwedge^* \mathfrak{p}$ v\'erifiant l'in\'egalit\'e
(\ref{eg}) est n\'ecessairement de la forme $\mu (\mathfrak{q} ')$ pour un certaine sous-alg\`ebre parabolique
$\theta$-stable $\mathfrak{q} ' \subset \mathfrak{g}$. Le $K$-type $\nu \subset \mathfrak{p} \wedge \mu (\subset
\bigwedge^{R+1} \mathfrak{p} )$ est donc \'egal \`a un certain $\mu (\mathfrak{q} ')$, o\`u $\mathfrak{q}'$ est une
sous-alg\`ebre parabolique $\theta$-stable de $\mathfrak{g}$.

\begin{lem} \label{l1}
Soient $\mathfrak{q}=\mathfrak{u} + \mathfrak{l}$ et $\mathfrak{q}'=\mathfrak{u}' +\mathfrak{l}'$ deux sous-alg\`ebres
paraboliques $\theta$-stable de $\mathfrak{g}$ telles que $L$ et $L'$ n'aient pas de facteurs simples compacts.
Supposons que le $K$-type $\mu (\mathfrak{q} ')$ intervienne dans
$\mathfrak{p} \wedge \mu (\mathfrak{q} )$. On est alors dans l'une des trois situations suivantes.
\begin{enumerate}
\item Le groupe $L$ a un centre non compact et $L'=L$.
\item Le groupe $L$ a un facteur simple localement isomorphe \`a $SO(n,1)$ ($n\geq 2$) ou \`a $SU(n,1)$ ($n\geq 1$) et
le groupe $L'$ est obtenu \`a partir de $L$ en rempla\c{c}ant ce facteur
par $SO(n-2,1) \times SO(2)$ ou $S(U(n-1,1) \times U(1))$ respectivement.
\item Le groupe $L'$ a un facteur simple localement isomorphe \`a $SO(n,1)$ ($n\geq 2$) ou \`a $SU(n,1)$ ($n\geq 1$)
et $L$ est obtenu \`a partir de $L'$ en rempla\c{c}ant ce facteur par $SO(n-2,1) \times SO(2)$ ou $S(U(n-1,1) \times U(1))$ respectivement.
\end{enumerate}
\end{lem}
{\it D\'emonstration du Lemme \ref{l1}.} Le $K$-type $\mu (\mathfrak{q} ') \subset \mathfrak{p} \wedge \mu (\mathfrak{q})$
intervient dans le produit tensoriel $\mathfrak{p} \otimes \mu (\mathfrak{q})$ son plus haut poids est donc de la
forme $2\rho (\mathfrak{u} \cap \mathfrak{p} ) + \delta$, o\`u $\delta$ est (un poids de la repr\'esentation de
$K$ dans $\mathfrak{p}$, c'est-\`a-dire) ou bien nulle ou bien une racine de
$\mathfrak{t}$ dans $\mathfrak{p}$. On a en fait seulement trois possibilit\'es~:
\begin{enumerate}
\item ou bien $\delta$ est nulle;
\item ou bien $\delta$ est une racine de $\mathfrak{t}$ dans $\mathfrak{l} \cap \mathfrak{p}$;
\item ou bien $-\delta$ est une racine de $\mathfrak{t}$ dans $\mathfrak{u} \cap \mathfrak{p}$.
\end{enumerate}

Nous allons montrer que ces trois cas correspondent aux diff\'erentes situations pr\'edites par le Lemme \ref{l1}.

Commen\c{c}ons par supposer que $\delta =0$. Alors (cf. \cite[Corollary 3.7]{VoganZuckerman}),
$${\rm Hom}_{L\cap K} (\mathfrak{l} \cap \mathfrak{p}, {\Bbb C} ) \cong {\rm Hom}_K (\bigwedge^{R+1} \mathfrak{p} ,
\mu (\mathfrak{q} )) \neq 0 ,$$
et le groupe $L$ a un centre non compact. Remarquons finalement que dans ce cas $\mu (\mathfrak{q}' ) = \mu (
\mathfrak{q} )$ et donc que $L=L'$. Nous sommes donc dans le premier cas pr\'edit par le Lemme \ref{l1}.

Supposons maintenant que $\delta$ est une racine de $\mathfrak{t}$ dans $\mathfrak{l} \cap \mathfrak{p}$.
Soit $V$ l'espace de dimension $1$~: $\bigwedge^R (\mathfrak{u} \cap \mathfrak{p} )$.
D'apr\`es \cite[Lemme 3.4]{VoganZuckerman}, un vecteur de poids $2\rho (\mathfrak{u} \cap \mathfrak{p} ) + \alpha$
dans $\bigwedge^{R+1} \mathfrak{p}$, avec $\alpha$ un poids de $\bigwedge^{*} \mathfrak{p}$ intervient n\'ecessairement
dans $V \otimes \bigwedge^* (\mathfrak{l} \cap \mathfrak{p})$. Le vecteur
de plus haut poids $\bigwedge^{R+1} (\mathfrak{u}' \cap \mathfrak{p})$ de $\mu (\mathfrak{q} ')$ est donc n\'ecessairement de la forme
$$\bigwedge {}^R (\mathfrak{u} \cap \mathfrak{p} ) \otimes v$$
o\`u $v$ est un vecteur de poids $\delta$ dans $\mathfrak{l} \cap \mathfrak{p}$. L'alg\`ebre $\mathfrak{u}$ est
donc contenue dans $\mathfrak{u}'$ et l'intersection $\mathfrak{q}'' = \mathfrak{q}' \cap \mathfrak{l}$ d\'efinit 
une sous-alg\`ebre parabolique $\theta$-stable de $\mathfrak{l}$ telle que
\begin{enumerate}
\item l'alg\`ebre $\mathfrak{u}'' \cap (\mathfrak{l} \cap \mathfrak{p})$ soit de dimension un, et
\item le groupe $L'=L''$.
\end{enumerate}
La premi\`ere condition est tr\`es restrictive, elle implique que le groupe $L$ a un facteur simple
localement isomorphe \`a $SO(n,1)$ ($n\geq 2$) ou \`a $SU(n,1)$ ($n\geq 1$). La seconde condition implique quant \`a elle
que le groupe $L'$ est obtenu \`a partir de $L$ en rempla\c{c}ant ce facteur
par $SO(n-2,1) \times SO(2)$ ou $S(U(n-1,1) \times U(1))$ respectivement.
Nous sommes donc dans le deuxi\`eme cas pr\'edit par le Lemme \ref{l1}.

Supposons finalement que $-\delta$ est une racine de $\mathfrak{t}$ dans $\mathfrak{u} \cap \mathfrak{p}$.
Alors $\mu (\mathfrak{q} ) = 2 \rho (\mathfrak{u}' \cap \mathfrak{p} ) + \alpha$ o\`u $\alpha$ est
une racine de $\mathfrak{t}$ dans $\mathfrak{l}' \cap \mathfrak{p}$. Nous sommes donc ramen\'e au cas pr\'ec\'edent d'o\`u
l'on d\'eduit que le groupe $L'$ a un facteur simple localement isomorphe \`a \`a $SO(n,1)$ ($n\geq 2$) ou \`a $SU(n,1)$
($n\geq 1$) et que le groupe $L$ est obtenu \`a partir de $L'$ en rempla\c{c}ant ce facteur
par $SO(n-2,1) \times SO(2)$ ou $S(U(n-1,1) \times U(1))$ respectivement. Nous sommes donc dans le troisi\`eme cas
pr\'edit par le Lemme \ref{l1}.

\bigskip

Concluons maintenant la d\'emonstration de l'implication inverse dans le Th\'eor\`eme \ref{rep isol} toujours sous l'hypoth\`ese que
les images $d({\rm Hom}_K (\mu , \pi_j ))$ sont non triviales.

On sait alors qu'il existe un $K$-type $\nu \subset \mathfrak{p} \wedge \mu (\subset
\bigwedge^{R+1} \mathfrak{p} )$ \'egal \`a un certain $\mu (\mathfrak{q} ')$, o\`u $\mathfrak{q}'$ est une
sous-alg\`ebre parabolique $\theta$-stable de $\mathfrak{g}$. Le Lemme \ref{l1} s'applique donc et on a trois
possibilit\'es. Si $L$ a un centre non compact, $L$ admet un facteur simple localement isomorphe au groupe $SO(1,1)$
ce qui contredit le premier point du Th\'eor\`eme \ref{rep isol}. Si $L$ a un facteur simple localement isomorphe \`a $SO(n,1)$ ($n\geq 2$) ou \`a $SU(n,1)$ ($n\geq 1$), c'est \'evidemment encore en contradiction avec le premier point du Th\'eor\`eme \ref{rep isol}. Si enfin il existe un groupe $L'$ ayant un facteur simple localement isomorphe \`a $SO(n,1)$
($n\geq 2$) ou \`a $SU(n,1)$ ($n\geq 1$) et tel que le groupe $L$ soit obtenu \`a partir de $L'$ en rempla\c{c}ant ce
facteur par $SO(n-2,1) \times SO(2)$ ou $S(U(n-1,1) \times U(1))$ respectivement, c'est \'evidemment en contradiction
avec le premier point du Th\'eor\`eme \ref{rep isol} si $n\geq 3$ dans le cas $SO$ et $n\geq 2$ dans le cas $SU$.
Puisque $SO(2,1)$ et $SU(1,1)$ sont tous les deux localement isomorphes au groupe $SL(2,{\Bbb R})$, il nous reste la
possibilit\'e que $L'$ admette un facteur simple localement isomorphe \`a $SL(2,{\Bbb R})$ et que le groupe $L$
soit obtenu \`a partir de $L'$ en rempla\c{c}ant ce facteur par son tore compact. Mais ce facteur de $L'$ correspond \`a
une racine imaginaire $\beta \in \Pi$ orthogonale aux racines dans $\Pi (\mathfrak{l})$ et le Th\'eor\`eme \ref{rep isol}
est d\'emontr\'e.

\bigskip

Supposons maintenant que l'image $d({\rm Hom}_K (\mu , \pi_j ))$ soit triviale \`a partir d'un certain rang $j\geq j_0$.
Supposons pour simplifier $j_0 =0$.

On munit chaque cha\^{\i}ne $C^* (\pi_j ) \subset \bigwedge^* \mathfrak{p}^* \otimes \pi_j$ du produit
scalaire obtenu en tensorisant le produit scalaire sur $\pi_j$ avec le produit scalaire sur $\bigwedge^* \mathfrak{p}^*$
d\'efini par la forme de Killing de $\mathfrak{g}$. Notons alors $\partial$ l'adjoint de l'op\'erateur $d$  pour ce
produit scalaire. Rappelons qu'il d\'ecoule du Lemme de Kuga, cf. \cite{BorelWallach}, que
si $\partial ({\rm Hom}_K (\mu , \pi_j )) =0$ alors (puisque $d({\rm Hom}_K (\mu , \pi_j ))=0$)
le casimir $\pi_j (\Omega )=0$. Mais dans ce cas et puisque $\mu = \mu (\mathfrak{q})$
intervient dans $\pi_j$, il d\'ecoule de \cite[Theorem 1.4]{VoganZuckerman} que la repr\'esentation $\pi_j$ est
n\'ecessairement \'egale \`a $\pi$ contrairement \`a notre choix de $\pi_j$. Les images
$\partial ({\rm Hom}_K (\mu , \pi_j ))$ sont donc toutes non triviales.

De mani\`ere duale \`a notre premier cas et toujours
quitte \`a extraire une sous-suite de $(\pi_j )$, on peut alors
supposer qu'il existe un $K$-type $\nu \subset \bigwedge^{R-1} \mathfrak{p}$ qui intervient dans
chaque $\pi_j$ et tel que $\mu \subset \mathfrak{p} \wedge \nu$.
Le Lemme \ref{dirac}
appliqu\'e \`a chacune des repr\'esentations $\pi_j$ et au $K$-type $\nu$ implique que le plus haut poids $\chi$
de $\nu$ v\'erifie l'\'egalit\'e
\begin{eqnarray} \label{eg}
||\rho ||^2 -||w(\chi - \rho_n ) + \rho_c ||^2 =0.
\end{eqnarray}
Ce qui, en utilisant l\`a encore \cite{Kumaresan}, montre que le $K$-type $\nu \subset
\bigwedge^{R-1} \mathfrak{p} $ est \'egal \`a un certain $\mu (\mathfrak{q} ')$, o\`u $\mathfrak{q}'$ est une
sous-alg\`ebre parabolique $\theta$-stable de $\mathfrak{g}$.

Le Lemme \ref{l1} s'applique donc \`a nouveau (en \'echangeant les r\^oles de $\mathfrak{q} $ et $\mathfrak{q} '$)
et on a trois
possibilit\'es. Soit $L=L'$ a un centre compact,
soit $L'$ a un facteur simple localement isomorphe \`a $SO(n,1)$ ($n\geq 2$) ou \`a $SU(n,1)$ ($n\geq 1$) et
le groupe $L$ est obtenu \`a partir de $L'$ en rempla\c{c}ant ce facteur par $SO(n-2,1) \times SO(2)$ ou $S(U(n-1,1) \times U(1))$ respectivement, soit le groupe $L$ a un facteur simple localement isomorphe \`a $SO(n,1)$
($n\geq 2$) ou \`a $SU(n,1)$ ($n\geq 1$). On retrouve donc les trois cas trait\'es ci-dessus qui m\`enent aux trois
m\^eme contradictions des hypoth\`eses du Th\'eor\`eme \ref{rep isol}. L'implication inverse est donc compl\`etement d\'emontr\'ee.

\bigskip

Montrons maintenant l'implication directe. Remarquons d'abord que la
d\'em- onstration de celle-ci par Vogan est tr\`es courte et imm\'ediate une fois admises les propri\'et\'es de l'induction
cohomologiques. Pour clarifier la suite du texte nous reprenons la d\'emonstration de Vogan en termes d'induction parabolique.
Commen\c{c}ons donc par quelques rappels sur l'induction parabolique. Soit $H$ un sous-groupe de Cartan $\theta$-invariant dans $G$. Alors 
$H=T^+ A$ (produit direct) avec $T^+ = H \cap K$ compact et $A= \exp (\mathfrak{h}_0 \cap \mathfrak{p}_0 )$. Nous notons 
$\Delta ( \mathfrak{g} , \mathfrak{h} )$ l'ensemble des racines de $\mathfrak{h}$ dans $\mathfrak{g}$. 

Soit ${\Bbb M}=MA$ la d\'ecomposition de Langlands du centralisateur de $A$ dans $G$. Le groupe $M$ est un groupe lin\'eaire r\'eductif et $T^+$ est un sous-groupe
de Cartan compact de $M$. Rappelons qu'alors $M$ poss\`ede une s\'erie discr\`ete. Plus pr\'ecisemment (cf. \cite{Knapp}), les repr\'esentations de 
la s\'erie discr\`ete de $M$ sont param\'etr\'ees par leur param\`etre d'Harish-Chandra $\lambda \in (\mathfrak{t}^+)^*$ associ\'e \`a un (pseudo-)caract\`ere r\'egulier de $T^+$ et 
regard\'e modulo l'action du groupe de Weyl $W(G/T^+ )$. Nous notons $\sigma (\lambda )$ la repr\'esentation correspondante.

Fixons $\gamma = (\lambda , \nu) \in \mathfrak{h}^*$ avec $\lambda$ comme au-dessus et $\nu \in \widehat{A}$ (nous notons \'egalement 
$\nu$ l'\'el\'ement de $\mathfrak{a}^*$ correspondant). On v\'erifie que ${\Bbb M}$ est le sous-groupe de Levi
d'un sous-groupe parabolique $P=MAN$ tel que 
$${\rm Re} \langle \alpha , \nu \rangle \leq 0 , \; \mbox{ pour toute racine } \alpha \in \Delta (\mathfrak{n} , \mathfrak{a} ).$$
La {\it repr\'esentation standard de param\`etre $\gamma$} est le $(\mathfrak{g}, K)$-module
$$X(\gamma ) = X^G (\gamma ) = {\rm ind}_P^G (\sigma (\lambda ) \otimes \nu \otimes 1) .$$
(Ici nous utilisons l'induction normalis\'ee et ne consid\'erons que l'espace des vecteurs $K$-finis de la repr\'esentation induite.)
Les diff\'erents constituants (sous-quotients) irr\'eductibles de $X(\gamma )$ ne d\'ependent pas de $P$
mais seulement de $\gamma$. Speh et Vogan associent dans \cite{SpehVogan} une famille canonique $\{
\overline{X}^1 (\gamma ) , \ldots , \overline{X}^r (\gamma ) \}$ de sous-quotients irr\'eductibles de $X(\gamma )$,
les {\it sous-quotients de Langlands}. Ce sont exactement les sous-quotients qui contiennent un $K$-type minimal de
$X(\gamma )$, il sera ainsi plus commode de les param\'etrer $\overline{X} (\gamma , \mu )$
par les diff\'erents $K$-types minimaux $\mu$ de $X(\gamma )$.

Cette construction est reli\'ee aux repr\'esentations cohomologiques $A_{\mathfrak{q}}$ de la mani\`ere suivante. Soit
$\mathfrak{q}= \mathfrak{l} + \mathfrak{u}$ une sous-alg\`ebre parabolique $\theta$-stable de $\mathfrak{g}$. Fixons $H=T^+A$ un sous-groupe de
Cartan $\theta$-stable et maximalement d\'eploy\'e dans $L$ et $L=(L\cap K ) A N^L$ une d\'ecomposition d'Iwasawa.
Nous voulons lui associer une repr\'esentation standard et donc un param\`etre $\gamma = (\lambda , \nu )$. On prend
$\nu \in \mathfrak{a}^*$ \'egal \`a la demi-somme des racines de $\mathfrak{a}$ dans $\mathfrak{n}^L$ et
$$\lambda = \rho^+ + \rho (\mathfrak{u} ) \in (\mathfrak{t}^+ )^* ,$$
o\`u un syst\`eme positif de racines de $\mathfrak{t}^+$ dans $\mathfrak{m} \cap \mathfrak{l}$ \'etant fix\'e, nous notons
$\rho^+$ la demi-somme des racines positives de $\mathfrak{t}^+$ dans $\mathfrak{m} \cap \mathfrak{l}$ et $\rho (\mathfrak{u})$
la demi-somme des racines de $\mathfrak{t}^+$ dans $\mathfrak{u}$. Alors, cf. \cite{VoganZuckerman}, la repr\'esentation $A_{\mathfrak{q}}$ s'identifie \`a
l'{\bf unique} ``quotient de Langlands'' $\overline{X} (\gamma )$ de la repr\'esentation induite $X(\gamma )$.

Le lien avec le sous-groupe $L$ est le suivant. Supposons que la sous-alg\`ebre de Cartan $\mathfrak{h}$ soit contenue dans la sous-alg\`ebre de Levi $\mathfrak{l}$
d'une sous-alg\`ebre parabolique $\theta$-stable $\mathfrak{q}=\mathfrak{l} + \mathfrak{u}$. Alors, l'ensemble des constituants de $X(\gamma )$ est en bijection
avec l'ensemble des constituants de $X^L (\gamma_L )$, o\`u $\gamma_L = (\lambda_L , \nu)$ avec $\lambda_L = \lambda - \rho (\mathfrak{u} )$. Si, plus pr\'ecisemment, $\mu_L - 2 \rho (\mathfrak{u} \cap \mathfrak{p} )$ est un
$(K \cap L)$-type minimal de $X^L (\gamma_L )$, il existe un unique $K$-type $\mu$ tel que $\mu_L$
soit la repr\'esentation de $L \cap K$ engendr\'ee par un espace de poids extr\`eme dans $\mu$. La
bijection associe alors \`a $\overline{X}^L (\gamma_L , \mu_L - 2\rho (\mathfrak{u} \cap \mathfrak{p} ))$
le module $\overline{X}^G (\gamma , \mu )$. Remarquons que ce lien
s'applique en particulier au cas des repr\'esentations cohomologiques. Dans ce cas $X^L (\gamma_L )$ a encore un unique ``quotient de Langlands''  $\overline{X}^L (\gamma_L )$ qui est la repr\'esentation triviale de $L$.
% de caract\`ere infinit\'esimal $(\rho^+ , \nu )= \rho_L$.

Revenons maintenant \`a la d\'emonstration de l'implication directe dans le Th\'eor\`eme \ref{rep isol} et
supposons d'abord que le groupe $L$ poss\`ede un facteur simple localement isomorphe \`a $SO(n,1)$ ou $SU(n,1)$ ($n\geq 1$). Alors la repr\'esentation triviale de $L$
n'est {\bf pas} isol\'ee dans le dual unitaire $\widehat L$ de $L$ ($L$ n'a pas la propri\'et\'e (T) de
Kazhdan)~: si $L$ poss\`ede un facteur simple localement isomorphe \`a $SO(1,1)$ il suffit de consid\'erer
une suite de caract\`eres unitaires convergeant vers le caract\`ere trivial, sinon $1_L$ est
limite de repr\'esentations de la s\'erie compl\'ementaire unitaire de $SO(n,1)$ ($n\geq 2$) ou $SU(n,1)$ ($n\geq 1$).
Dans tous les cas la repr\'esentation triviale de $L$ est approch\'ee par des modules $\overline{X}^L (\gamma_L^i ,
\mu_L^i)$ auxquels sont naturellement associ\'es les modules $\overline{X} (\gamma^i , \mu^i)$ qui sont unitaires
d'apr\`es \cite{Vogan2} et convergent vers la repr\'esentation cohomologique $A_{\mathfrak{q}}$.
%ce sont les nu qui doivent bouger (caracteres !)

Supposons enfin qu'il y ait une racine imaginaire non compacte $\beta$ dans $\Pi$ orthogonale \`a $\Pi (\mathfrak{l} )$.
Consid\'erons alors la sous-alg\`ebre parabolique $\theta$-stable $\mathfrak{q}' = \mathfrak{l}' + \mathfrak{u}'$
correspondant \`a $\Pi (\mathfrak{l} ) \cup \{ \beta \}$. Le sous-groupe de Levi $L'$ est localement isomorphe au
produit $L \times SL(2,{\Bbb R})$ au centre pr\`es. La correspondance
ci-dessus associe \`a la repr\'esentation triviale $\overline{X}^L (\gamma_L )$ de $L$ une
repr\'esentation $\overline{X}^{L'} (\gamma_{L'} )$ de $L'$, dont on peut facilement v\'erifier que c'est
la premi\`ere s\'erie discr\`ete de $SL(2,{\Bbb R})$ (tensoris\'ee avec la repr\'esentation triviale de $L$), \`a
laquelle correspond, via la correspondance de $L'$ \`a $G$ et par fonctorialit\'e de l'induction,
la repr\'esentation cohomologique $A_{\mathfrak{q}}$.
Mais la premi\`ere s\'erie discr\`ete de $SL(2,{\Bbb R})$ est dans l'adh\'erence de la s\'erie compl\'ementaire et comme
au paragraphe pr\'ec\'edent la correspondance de $L'$ \`a $G$ permet de construire une suite de modules unitaires qui convergent vers $A_{\mathfrak{q}}$. Ceci conclut la d\'emonstration du Th\'eor\`eme \ref{rep isol}.

\bigskip

\noindent
{\bf Remarques.} {\bf 1.} Dans \cite[Theorem A.10]{Vogan}, Vogan \'enonce en fait quatre conditions n\'ecessaires et suffisantes pour que la
repr\'esentation $\pi$ du Th\'eor\`eme \ref{rep isol} soit isol\'ee. La condition 0) de Vogan est que le groupe
$L$ n'ait pas de facteurs simples compacts (non ab\'elien), dans le cas des repr\'esentations cohomologiques
auxquelles nous nous sommes restreint ici, nous avons rappel\'e (Th\'eor\`eme \ref{rep cohom}) que l'on pouvait
toujours supposer cette condition v\'erifi\'ee. La condition 1) de Vogan est que le centre de $L$ ne soit pas compact,
comme le remarque Vogan cette condition \'equivaut \`a ce que le groupe $L$ n'ait pas de facteur localement isomorphe
au groupe $SO(1,1)$, elle est donc contenue dans notre condition 1. La condition 2) de Vogan correspond \`a notre
condition 1. (avec $n\geq 2$ dans le cas $SO$). Enfin la condition 3) de Vogan peut para\^{\i}tre plus faible puisque
dans notre cas celle-ci s'\'ecrit~:
\begin{center} {\it
$$\langle \beta^{\vee} , \rho_{\mathfrak{g}} \rangle \neq 1$$
pour toute racine (imaginaire) non compacte $\beta \in \Pi$ orthogonale \`a $\Pi (\mathfrak{l} )$.}
\end{center}
Cette condition est en fait \'equivalente \`a notre condition 2. En effet, si $\beta$ est une racine simple
la transformation $s_{\beta}$ permute les racines positives $\neq \beta$ et envoie $\beta$ sur $-\beta$. Elle envoie
donc la demi-somme $\rho_{\mathfrak{g}}$ de toutes les racines positives sur $\rho_{\mathfrak{g}} - \beta$. Puisque
par ailleurs $s_{\beta} ( \alpha ) = \alpha - 2 \langle \beta^{\vee} , \alpha \rangle \beta$, on obtient que
pour toute racine $\beta \in \Pi$, $\langle\beta^{\vee} , \rho_{\mathfrak{g}} \rangle = 1$.

{\bf 2.} Si dans le Th\'eor\`eme \ref{rep isol}, la repr\'esentation $\pi$ appartient \`a la s\'erie discr\`ete de $G$,
elle n'est pas isol\'ee. En effet, l'alg\`ebre $\mathfrak{l}$ est alors compacte et donc ab\'elienne  et toute
racine $\beta \in \Pi$ est orthogonale \`a $\Pi (\mathfrak{l} )$.

{\bf 3.} Dans la d\'emonstration qu'il donne du Th\'eor\`eme \ref{rep isol}, Vogan distingue deux cas suivant que le
$K$-type $\mu$ soit ou non
un $K$-type minimal de toutes (sauf peut-\^etre un nombre fini d'entre elles) les repr\'esentations $\pi_j$.
La d\'emonstration est plus simple lorsque $\mu$ est effectivement un $K$-type minimal des
repr\'esentations $\pi_j$. Et sous cette hypoth\`ese la repr\'esentation $\pi$ est isol\'ee si et seulement si
la condition 1. du Th\'eor\`eme \ref{rep isol} est v\'erifi\'ee. Ceci correspond au fait que dans notre Lemme
\ref{l1}, pour que $\mu (\mathfrak{q}' )$ ne soit pas un plus bas $K$-type que $\mu (\mathfrak{q} )$, l'on doive
exclure le troisi\`eme cas.

{\bf 4.} Pr\'ecisons la correspondance entre certains $(\mathfrak{l} , L \cap K)$-modules et $(\mathfrak{g} , K)$-modules utilis\'ee dans la d\'emonstration de
l'application directe du Th\'eor\`eme \ref{rep isol}. Soit $\mathfrak{q} = \mathfrak{l} + \mathfrak{u}$ une sous-alg\`ebre parabolique $\theta$-stable de $\mathfrak{g}$.
Fixons une sous-alg\`ebre de Cartan $\mathfrak{h} \subset \mathfrak{l}$, et un poids $\gamma_L \in \mathfrak{h}^*$. Remarquons que 
$$\rho (\mathfrak{u} ) = \frac12 \sum_{\alpha \in \Delta (\mathfrak{g} , \mathfrak{h} )} \alpha \in \mathfrak{h}^* .$$
Soit $Y$ un $(\mathfrak{l} , L \cap K)$-module irr\'eductible de caract\`ere infinit\'esimal $\gamma_l  = \gamma - \rho (\mathfrak{u} )$. D'apr\`es Langlands \cite{Langlands} et 
Knapp-Zuckerman \cite{KnappZuckerman}, il existe $\lambda \in (\mathfrak{t}^+)^*$ associ\'e \`a un (pseudo-)caract\`ere r\'egulier de $T^+$ et $\nu \in \widehat{A}$ tels que
$\gamma = (\lambda , \nu ) \in \mathfrak{h}^*$ et $Y$ soit l'unique quotient de Langlands $\overline{X}^L (\gamma_L )$ de $X^L (\gamma_L )$.  Supposons de plus 
\begin{eqnarray} \label{conda}
{\rm Re} \langle \alpha , \gamma \rangle \geq 0 , \; \mbox{ pour tout } \alpha \in \Delta (\mathfrak{u} , \mathfrak{h} ),
\end{eqnarray}
et 
\begin{eqnarray} \label{condb}
\langle \alpha , \gamma \rangle \neq 0 , \; \mbox{ pour tout } \alpha \in \Delta (\mathfrak{u} , \mathfrak{h} ).
\end{eqnarray}
Alors le module $\overline{X} (\gamma )$ est bien d\'efinit, on le note ${\cal R} Y$. C'est un $(\mathfrak{g} , K)$-module irr\'eductible de caract\`ere infinit\'esimal
$\gamma$ naturellement attach\'e \`a $Y$. Le foncteur ${\cal R}$ est
appel\'e {\it foncteur d'induction cohomologique}, il r\'ealise une equivalence de cat\'egorie exacte
sur son image. Il envoie repr\'esentations unitaires sur repr\'esentations unitaires et repr\'esentations non-unitaires sur 
repr\'esentations non-unitaires, cf. \cite{Vogan2}. 

{\bf 5.} La difficult\'e de la d\'emonstration originale de Vogan tient \`a ce qu'une suite $(\pi_j ) \subset \widehat{G}$ adh\'erente \`a ${\cal R} 1_L$ n'appartient pas n\'ecessairement \`a l'image de ${\cal R}$.
La d\'emonstration de Vogan implique n\'eanmoins qu'il existe alors $L'$ li\'e \`a $L$ comme dans le point 3. du Lemme \ref{l1} 
avec (en particulier $L'$ poss\'edant un facteur de rang $1$ et donc) $1_{L'}$ non-isol\'ee~: $\rho_i \rightarrow 1_{L'}$ telle
que la suite $({\cal R}' \rho_i )_i$ adh\`ere \`a ${\cal R} 1_L$. O\`u nous avons not\'e ${\cal R}'$ le foncteur d'induction cohomologique de $L'$ \`a $G$.

\bigskip

Dans \cite{Livre} et \cite{Lefschetz} une propri\'et\'e d'isolation faible joue un r\^ole important~:
l'isolation sous la condition $d=0$. Nous dirons d'une repr\'esentation cohomologique $\pi$ comme dans le Th\'eor\`eme
\ref{rep isol} qu'elle est {\it isol\'ee sous la condition $d=0$} si elle est isol\'ee de l'ensemble des repr\'esentations
irr\'eductibles unitaires telles que ${\rm Im} (d_R )=0$, o\`u $d_R$
d\'esigne la diff\'erentielle sur les cocha\^{\i}nes de
degr\'e $R= R(\mathfrak{q} )$ du complexe calculant la $(\mathfrak{g} , K)$-cohomologie.

Notre d\'emonstration du Th\'eor\`eme \ref{rep isol} est bien adapt\'ee pour caract\'eriser les repr\'esentations
cohomologiques isol\'ees sous la condition $d=0$. On obtient plus pr\'ecisemment.

\begin{prop} \label{d}
Soit $\pi$ comme dans le Th\'eor\`eme \ref{rep cohom}. La repr\'esentation $\pi$ est {\bf isol\'ee sous la condition
$d=0$} si et seulement s'il n'existe aucune sous-alg\`ebre parabolique $\theta$-stable $\mathfrak{q}' = \mathfrak{u}'
+ \mathfrak{l}'$ telle que le groupe $L'$ ait un facteur simple localement isomorphe \`a $SO(n,1)$ ($n\geq 2$) ou
\`a $SU(n,1)$ ($n\geq 1$) et que le groupe $L$ s'obtienne \`a partir de $L'$ en rempla\c{c}ant ce facteur par $SO(n-2,1)
\times SO(2)$ ou $S(U(n-1,1) \times U(1))$ respectivement.
\end{prop}
{\it D\'emonstration.} Montrons que c'est une condition suffisante par contrapos\'ee. Supposons donc qu'il existe
une suite $(\pi_j)$ comme dans la d\'emonstration du Th\'eor\`eme \ref{rep isol} et v\'erifiant $d(C^R (\pi_j ))=0$.
Alors $d({\rm Hom}_K (\mu , \pi_j ))=0$ et la d\'emonstration du Th\'eor\`eme \ref{rep isol} implique qu'il existe
un $K$-type $\mu ' = \mu (\mathfrak{q}')$ associ\'e \`a une certaine sous-alg\`ebre parabolique $\theta$-stable de
$\mathfrak{g}$ tel que $\mu ' \subset \bigwedge^{R-1} \mathfrak{p}$ intervienne dans
chaque $\pi_j$ et tel que $\mu \subset \mathfrak{p} \wedge \mu '$. Le Lemme \ref{l1} s'applique en \'echangeant
les r\^oles de $\mathfrak{q}$ et $\mathfrak{q}'$. En outre, puisque $\mu ' \subset \bigwedge^{R-1} \mathfrak{p}$,
le plus haut poids $2\rho (\mathfrak{u} \cap  \mathfrak{p} )$ de $\mu$ est n\'ecessairement de la forme
$2\rho (\mathfrak{u} ' \cap \mathfrak{p} ) + \delta$ avec $\delta$ racine {\bf positive}. On est donc dans le deuxi\`eme
cas du Lemme \ref{l1} (en \'echangeant les r\^oles de $L$ et $L'$). Ce qui implique la Proposition \ref{d}.

\bigskip

\section{Illustration dans le cas des groupes classiques}

Dans cette section, nous appliquons le Th\'eor\`eme \ref{rep isol} et la Proposition \ref{d} aux
diff\'erents groupes classiques.

\paragraph{Groupes unitaires.} Soient $p$ et $q$ des entiers strictement positifs avec $p\leq q$ et
\begin{eqnarray} \label{Gnc}
G = U(p,q) := \left\{ g= \left(
\begin{array}{cc}
A & B \\
C & D
\end{array} \right) \; : \; ^t \! \overline{g} \left(
\begin{array}{cc}
1_p & 0 \\
0 & -1_q
\end{array} \right) g = \left(
\begin{array}{cc}
1_p & 0 \\
0 & -1_q
\end{array} \right) \right\} ,
\end{eqnarray}
o\`u $A\in M_{p\times p} ({\Bbb C})$, $B\in M_{p\times q} ({\Bbb C})$, $C\in M_{q\times p} ({\Bbb C})$ et
$D\in M_{q\times q} ({\Bbb C})$. Le rang r\'eel de $G$ est alors $p$ et le sous-groupe  
$$K = \left\{ g =  \left(
\begin{array}{cc}
A & 0 \\
0 & D
\end{array} \right) \in G \; : \;  A\in U(p) , \ D\in U(q) \right\} $$
est un sous-groupe compact maximal dans $G$.

Rappelons que la multiplication par $i =\sqrt{-1}$ induit une d\'ecomposition
$$\mathfrak{p} = \mathfrak{p}^+ \oplus \mathfrak{p}^- .$$
L'alg\`ebre de Lie $\mathfrak{g}$ est bien \'evidemment $M_{(p+q) \times (p+q)} ({\Bbb C})$, et l'on voit ses
\'el\'ements sous forme de blocs comme dans (\ref{Gnc}). On a alors,
$$\mathfrak{p}^+ = \left\{ \left(
\begin{array}{cc}
0 & B \\
0 & 0
\end{array} \right) \mbox{ avec } B \in M_{p\times q} ({\Bbb C} ) \right\}$$
et
$$\mathfrak{p}^- = \left\{ \left(
\begin{array}{cc}
0 & 0 \\
C & 0
\end{array} \right) \mbox{ avec } C \in M_{q \times p} ({\Bbb C} )\right\} .$$
Soit $E= {\Bbb C}^p$ (resp. $F= {\Bbb C}^q$) la repr\'esentation standard de $U(p)$ (resp. $U(q)$). Alors, comme
repr\'esentation de $K_{{\Bbb C}}$, $\mathfrak{p}^+ = E \otimes F^*$.

Dans \cite{Tentative} nous avons param\'etr\'e les repr\'esentations cohomologiques de $G$ par certains 
couple de partitions. 
Rappelons qu'une {\it partition} est une suite d\'ecroissante $\lambda$ d'entiers naturels $\lambda_1 \geq \ldots \geq
\lambda_l \geq 0$. Les entiers $\lambda_1 , \ldots , \lambda_l$ sont des {\it parts}. La {\it longueur} $l(\lambda )$ d\'esigne le
nombre de parts non nulles, et le {\it poids} $|\lambda |$, la somme des parts. On se soucie peu, d'ordinaire, des parts nulles~: on se
permet en particulier, le cas \'ech\'eant, d'en rajouter ou d'en \^oter.
Le {\it diagramme de Young} de $\lambda$, que l'on notera \'egalement $\lambda$, s'obtient en superposant, de haut en bas, des lignes
dont l'extr\'emit\'e gauche est sur une m\^eme colonne, et de longueurs donn\'ees par les parts de $\lambda$. Par sym\'etrie
diagonale, on obtient le diagramme de Young de la {\it partition conjugu\'ee}, que l'on notera $\lambda^*$.

Soient $\lambda$ et $\mu$ deux partitions telles que le diagramme de $\mu$ contienne $\lambda$, ce que nous noterons $\lambda \subset \mu$. Notons
$\mu / \lambda$ le compl\'ementaire du diagramme de $\lambda$ dans celui de $\mu$~: c'est une {\it partition gauche} son diagramme
est un {\it diagramme gauche}. Dans la pratique les partitions $\lambda$ que nous rencontrerons seront incluses dans la {\it partition
rectangulaire} $p\times q = (\underbrace{q, \ldots , q}_{p \; {\rm fois}}) =(q^p)$, le diagramme gauche $p\times q /\lambda$ est alors
le diagramme de Young d'une partition auquel on a appliqu\'e une rotation d'angle $\pi$; nous noterons $\hat{\lambda}$ cette partition,
la {\it partition compl\'ementaire} de $\lambda$ dans $p \times q$.

Nous dirons d'un couple de partitions $(\lambda , \mu )$ qu'il est {\it compatible} ({\it compatible dans $p\times q$}
en cas d'ambiguit\'e) si
\begin{enumerate}
\item on a la suite d'inclusion $\lambda \subset \mu \subset p\times q$, et
\item le diagramme gauche $\mu / \lambda$ est une r\'eunion de diagrammes rectangulaires $p_i \times q_i$, $i=1, \ldots ,m$
ne s'intersectant qu'en des sommets.
\end{enumerate}

Les r\'esultats de Parthasarathy, Kumaresan et Vogan-Zuckerman mentionn\'es plus haut affirment alors que l'ensemble des 
classes d'\'equivalences de repr\'esenta- tions cohomologiques de $G$ est 
$$\{ A(\lambda , \mu ) \; : \; (\lambda , \mu ) \mbox{ est un couple compatible de partitions} \},$$
o\`u 
$$H^k (\mathfrak{g} , K , A(\lambda , \mu )) \cong \left\{ 
\begin{array}{ll}
{\Bbb C} & \mbox{ si } k=R:= |\lambda| + |\hat{\mu}| = pq - \sum_{i=1}^m p_i q_i , \\ 
0 & \mbox{ si } k<R , 
\end{array} \right.$$
et plus g\'en\'eralement
$$H^k (\mathfrak{g} , K , A(\lambda , \mu )) \cong H^{k-R} \left( \prod_{i=1}^m U(p_i + q_i )/(U(p_i ) \times U(q_i )), {\Bbb C} \right).$$
Ici $U(p_i + q_i )/(U(p_i ) \times U(q_i ))$ est la grassmannienne des $p_i$-plans dans ${\Bbb C}^{p_i + q_i}$.

La repr\'esentation triviale correspond au couple $(0 , p\times q)$ et le degr\'e $R$ correspondant est \'evidemment \'egal 
\`a $0$. On retrouve que la cohomologie qui lui correspond est la cohomologie du dual compact 
$\widehat{X}_G = U(p + q )/(U(p ) \times U(q ))$ de $X_G$. 
Les repr\'esentations cohomologiques appartenant \`a la s\'erie discr\`ete sont celles qui v\'erifient $\mu = \lambda$. Enfin, remarquons
que le ``type de Hodge'' de la classe de cohomologie fortement primitive de $A(\lambda , \mu )$ est $(R^+, R^- ) = (|\lambda |, 
|\hat{\mu} | )$. Nous noterons $H^{\lambda , \mu}$ la $A(\lambda , \mu )$-composante de la cohomologie de degr\'e $R$ (composante fortement primitive).

Le Th\'eor\`eme \ref{rep isol} implique qu'un certain nombre de ces repr\'esentations sont isol\'ees dans le dual unitaire de $G$.
Il implique plus pr\'ecisemment le corollaire suivant.

\begin{cor} \label{isolU}
Soit $(\lambda , \mu)$ un couple compatible de partitions dans $p\times q$ avec $\mu / \lambda = (p_1 \times q_1 ) * \ldots *
(p_m \times q_m)$. Alors la repr\'esentation $A(\lambda , \mu )$ est isol\'ee dans le dual unitaire de $G$ si et seulement s'il n'existe
aucun couple compatible de partitions $(\lambda ' , \mu ')$ dans $p \times q$ tel que les diagrammes gauches $\mu ' /\lambda '$ et $\mu / \lambda$ ne diff\`erent
que d'une case. Elle est isol\'ee sous la condition $d=0$ si et seulement s'il n'existe aucun couple compatible de partitions $(\lambda ' , \mu ')$ dans $p\times q$ tel que 
le diagramme gauche $\mu / \lambda $ s'obtienne \`a partir de $\mu ' / \lambda ' $ en lui enlevant une case.
\end{cor}

On peut expliciter ces conditions~: la repr\'esentation $A(\lambda , \mu )$ est isol\'ee dans le dual unitaire de $G$ si et seulement si 
\begin{enumerate}
\item $\min_i (p_i , q_i ) \geq 2$, et 
\item si $\lambda_i = \mu_i > \lambda_{i+1}$ ($i=1 , \ldots , p$) alors $\mu_{i+1} = \mu_i$ (o\`u nous adoptons exceptionnellement la convention que 
$\lambda_{p+1} = \mu_{p+1} = -1$).
\end{enumerate}
Le dernier point signifie que $\lambda$ et $\mu$ n'ont aucun angle
$\begin{array}{c|}
\\ \hline
\end{array}$ ou 
$\begin{array}{|c}
\hline
\\ 
\end{array}$ en commun. En particulier, les repr\'esentations $A(\lambda , \mu)$ telles que $\lambda = \mu$ 
(qui sont exactement les repr\'esentations cohomologique de la s\'erie discr\`ete) ne sont {\bf jamais isol\'ees}.

\paragraph{Groupes orthogonaux.} Dans ce paragraphe $G = O(p,q)^0$, o\`u 
$p$ et $q$ sont des entiers strictement positifs. Le rang r\'eel de $G$ est donc $\min (p,q)$. On a 
\begin{eqnarray}
G = \left\{ g = \left(
\begin{array}{cc}
A & B \\
C & D 
\end{array} \right) \; : \; {}^t g \left(
\begin{array}{cc}
1_p & 0 \\
0 & -1_q
\end{array} \right) g = \left(
\begin{array}{cc}
1_p & 0 \\
0 & -1_q 
\end{array} \right) \right\}^0 ,
\end{eqnarray}
o\`u $A \in M_{p\times p} ({\Bbb R})$, $B \in M_{p\times q} ({\Bbb R})$, $C\in M_{q\times p} ({\Bbb R})$ et $D \in M_{q\times q } ({\Bbb R})$. Et, 
$$K = \left\{ g = \left( 
\begin{array}{cc}
A & 0 \\
0 & D 
\end{array} \right) \in G \; : \; A \in SO(p) , D \in SO(q) \right\} .$$
L'involution de Cartan $\theta$ est donn\'ee par $x \mapsto - {}^t x$. On a alors, 
$$\mathfrak{p} = \left\{ \left( 
\begin{array}{cc}
0 & B \\
{}^t B & 0
\end{array} \right) \; : \; B \in M_{p\times q} ({\Bbb C}) \right\} .$$
Soit $E= {\Bbb C}^p$ (resp. $F={\Bbb C}^q$) la repr\'esentation standard de $SO(p,{\Bbb C})$ (resp. $SO(q,{\Bbb C})$). Alors, comme repr\'esentation de 
$K_{{\Bbb C}}$, $\mathfrak{p} = E \otimes F^*$. 

Dans \cite{Lefschetz}, et toujours \`a partir des r\'esultats de Parthasarathy, Kumaresan et Vogan-Zuckerman, nous montrons que 
l'ensemble des 
classes d'\'equivalences de repr\'esentations cohomologiques de $G$ est alors
$$\{ A(\lambda )^{\pm_1}_{\pm_2} \; : \; (\lambda , \hat{\lambda} ) \mbox{ est un couple compatible de partitions} \footnote{Nous dirons
dor\'enavant d'une partition $\lambda$ telle que le couple $(\lambda , \hat{\lambda} )$ soit compatible, qu'elle est 
{\it orthogonale}.} \},$$
o\`u les signes $\pm_1$ et $\pm_2$ correspondent au fait que pour nous une classe d'\'equivalen- ce de repr\'esentations cohomologiques
correspond \`a une repr\'esentation du groupe connexe $G$ et que certaines repr\'esentations $A(\lambda )$ du groupe $O(p,q)$ 
se restreignent au groupe $G$ en une somme de $2$ ou $4$ repr\'esentations (cohomologiques) irr\'eductibles. On peut oublier
tout ceci dans le reste du texte.

Remarquons qu'une partition orthogonale $\lambda$ v\'erifie que son diagramme gauche associ\'e $\hat{\lambda} / \lambda$
s'\'ecrit comme une r\'eunion de diagrammes rectangulaires~: 
$(a_1 \times b_1) * \ldots *(a_m \times b_m ) * (p_0 \times q_0 )* (a_m \times b_m ) * \ldots * (a_1 \times b_1)$. 
Les groupes de $(\mathfrak{g} , K)$-cohomologie de $A(\lambda )^{\pm_1}_{\pm_2}$ se calculent alors de la mani\`ere suivante~: 
$$H^k (\mathfrak{g} , K , A(\lambda )^{\pm_1}_{\pm_2} ) 
 \cong \left\{ 
\begin{array}{ll}
{\Bbb C} & \mbox{ si } k=R:= |\lambda| , \\ 
0 & \mbox{ si } k<R , 
\end{array} \right. $$
et plus g\'en\'eralement
$$\begin{array}{l}
H^k (\mathfrak{g} , K , A(\lambda )^{\pm_1}_{\pm_2} ) \\
\cong H^{k-R} \left( O(p_0 + q_0)/(O(p_0 ) \times O(q_0 )) \prod_{i=1}^m U(a_i + b_i )/(U(a_i ) \times U(b_i )), {\Bbb C} \right).
\end{array}$$
Nous noterons $H^{\lambda}$ la $A(\lambda )$-composante de la cohomologie de degr\'e $R$.

Dans ce cas le Th\'eor\`eme \ref{rep isol} implique le corollaire suivant.

\begin{cor} \label{Oisol}
Soit $\lambda \subset p\times q$ une partition orthogonale avec $\hat{\lambda} / \lambda =
(a_1 \times b_1 ) * \ldots * (a_m \times b_m ) * (p_0 \times q_0 ) * (a_m \times b_m ) * \ldots *(a_1 \times b_1)$. Alors
les diff\'erentes repr\'esentations $A(\lambda )_{\pm_1}^{\pm_2}$ associ\'ees \`a $\lambda$ sont
soit toutes isol\'ees dans le dual unitaire de $O(p,q)^0$ soit toutes non isol\'ees. Elles sont isol\'ees si et seulement s'il n'existe aucune partition orthogonale 
$\lambda ' \subset p\times q$ telle que les diagrammes gauches $\hat{\lambda}/\lambda$ et $\widehat{\lambda '} / \lambda '$ ne diff\`erent que de deux cases.
Elle est isol\'ee sous la condition $d=0$ si et seulement s'il existe une partition orthogonale $\lambda ' \subset p\times q$ telle que
le diagramme gauche $\hat{\lambda}/\lambda$ s'obtienne \`a partir de $\widehat{\lambda '} / \lambda '$ en lui enlevant deux cases.
\end{cor}

\paragraph{Groupes $\mathbf{G=Sp(p,q)}$.} Ce cas est similaire \`a celui des groupes unitaires. Rappelons juste qu'en utilisant l'identification
classique entre l'alg\`ebre des quaternions ${\Bbb H}$ et ${\Bbb C}^2$, on a~:
$$\mathfrak{k}_0 = \left\{ \left(
\begin{array}{cccc}
A & 0 & S & 0 \\
0 & B & 0 & T \\
-\overline{S} & 0 & \overline{A} & 0 \\
0 & -\overline{T} & 0 & \overline{B} 
\end{array} \right) \; : \; 
\begin{array}{l}
A \in \mathfrak{u} (p) , \; B \in \mathfrak{u} (q), \\
S \in {\rm Sym}(p,{\Bbb C}) , \; T \in {\rm Sym} (q,{\Bbb C})
\end{array} \right\} , $$
$$\mathfrak{p}_0 = \left\{ \left(
\begin{array}{cccc}
0 & M & 0 & X \\
{}^t \overline{M} & 0 & {}^t X & 0 \\
0 & \overline{X} & 0 & -\overline{M} \\
{}^t \overline{X} & 0 & -{}^t M & 0 
\end{array} \right) \; : \; M ,X \in M_{p \times q} ({\Bbb C}) \right\} .$$ 
Comme dans le cas unitaire \cite{Tentative} il n'est alors pas difficile de montrer que l'ensemble des repr\'esentations cohomologiques de $G$ est 
$$\left\{ A(\lambda , \mu )_i \; : \;
\begin{array}{l}
(\lambda , \mu ) \mbox{ est un couple compatible de partitions et } \\
 i=0,1 \; (\mbox{toujours } =1 \mbox{ si } \lambda_p >0 \mbox{ ou } \mu_p = 0)
 \end{array}  \right\},$$
o\`u 
$$H^k (\mathfrak{g} , K , A(\lambda , \mu)_1 ) \cong \left\{ 
\begin{array}{ll}
{\Bbb C} & \mbox{ si } k=R:= pq +|\lambda| + |\hat{\mu}| = 2pq - \sum_{i=1}^m p_i q_i , \\ 
0 & \mbox{ si } k<R , 
\end{array} \right.$$
et 
$$H^k (\mathfrak{g} , K , A(\lambda , \mu)_0 ) \cong \left\{
\begin{array}{ll}
{\Bbb C} & \mbox{ si } k=R:= \left\{
\begin{array}{l}
pq - p_1 q_1 + |\lambda| + |\hat{\mu}| \\
2pq -2p_1 q_1- \sum_{i=2}^m p_i q_i ,
\end{array} \right. \\ 
0 & \mbox{ si } k<R , 
\end{array} \right.$$
et plus g\'en\'eralement
$$H^k (\mathfrak{g} , K , A(\lambda , \mu)_1 ) \cong H^{k-R} \left( \prod_{i=1}^m U(p_i + q_i )/(U(p_i ) \times U(q_i )), {\Bbb C} \right),$$
et
$$\begin{array}{l}
H^k (\mathfrak{g} , K , A(\lambda , \mu)_0 ) \\
\cong H^{k-R} \left( Sp (p_1 + q_1 ) /(Sp (p_1 ) \times Sp (q_1 )) \times \prod_{i=2}^m U(p_i + q_i )/(U(p_i ) \times U(q_i )), {\Bbb C} \right).
\end{array}$$

Dans ce cas le Th\'eor\`eme \ref{rep isol} implique le corollaire suivant.

\begin{cor} \label{isolSppq}
Soit $(\lambda , \mu)$ un couple compatible de partitions dans $p\times q$ avec $\mu / \lambda = (p_1 \times q_1 ) * \ldots *
(p_m \times q_m)$. Alors la repr\'esentation $A(\lambda , \mu )_i$ ($i=0,1$) 
est isol\'ee dans le dual unitaire de $G$ si et seulement si
\begin{enumerate}
\item lorsque $i=1$~: il n'existe aucun couple compatible de partitions $(\lambda ' , \mu ')$ dans $p\times q$ tel que 
les diagrammes gauches $\mu' / \lambda '$ et $\mu /\lambda$ ne diff\`erent que d'une case, et
\item lorsque $i=0$~:  $p_1 , q_1 \geq 1$, $p_1 + q_1 \geq 3$ et il n'existe aucun couple compatible de partitions $(\lambda ' , \mu ')$ dans $p\times q$ tel que 
les diagrammes gauches $\mu' / \lambda '= (p_1 ' \times q_1 ') * \ldots *(p_m ' \times q_m ')$ et $\mu /\lambda$ ne diff\`erent que d'une case et 
$(p_1 ' , q_1 ')=(p_1 , q_1)$. 
\end{enumerate}
Elle est isol\'ee sous la condition $d=0$ si et seulement s'il n'existe aucun couple compatible de partitions $(\lambda ' , \mu ')$ dans $p\times q$ tel que 
le diagramme gauche $\mu /\lambda$ s'obtienne \`a partir de $\mu ' / \lambda '$ en lui enlevant une case et $(p_1 , q_1) = (p_1 ' , q_1 ' )$ si $i=0$.
\end{cor}

\section{Isolation automorphe}

Consid\'erons \`a nouveau un groupe $G$ semi-simple alg\'ebrique sur ${\Bbb Q}$. 
Si l'on se restreint au dual automorphe il est naturel d'esp\'erer des propri\'et\'es d'isolation plus fortes. Nous conjecturons le r\'esultat suivant.

\begin{conj} \label{isolaut}
Soit $\pi$ comme dans le Th\'eor\`eme \ref{rep cohom}. La repr\'esentation $\pi$ est {\bf isol\'ee} dans 
$$\{ \pi \} \cup \widehat{G}_{{\rm Aut}}$$
d\`es que le sous-groupe de Levi $L \subset G$, associ\'e \`a la sous-alg\`ebre
parabolique $\mathfrak{q}$, a un centre compact. C'est en particulier toujours le cas lorsque ${\rm rang}_{\Bbb C} (G)
= {\rm rang}_{\Bbb C} (K)$, autrement dit lorsque $G$ poss\`ede une s\'erie discr\`ete. 
\end{conj}

Nous conjecturons de plus que la Conjecture \ref{isolaut} est optimale au sens faible suivant~:
le type r\'eel de $G$ \'etant fix\'e, il devrait exister une ${\Bbb Q}$-forme de $G$ telle que la repr\'esentation cohomologique $\pi$ soit
contenue et non-isol\'ee dans $\widehat{G}_{{\rm Aut}}$.

Enfin nous conjecturons qu'une repr\'esentation cohomologique $\pi$ de $G$ est toujours {\bf isol\'ee sous la condition $d=0$} dans 
$$\{ \pi \} \cup \widehat{G}_{{\rm Aut}}.$$

\medskip

Dans la suite de cette section nous cherchons \`a motiver la Conjecture \ref{isolaut} \`a l'aide des Conjectures d'Arthur. 
Dans deux articles fondamentaux, Arthur a en effet donn\'e une description conjecturale
des repr\'esentations des groupes r\'eductifs qui peuvent appara\^{\i}tre dans $L^2 (\Gamma \backslash
G)$ pour un sous-groupe de congruence. Avec Clozel dans \cite{Livre}, nous avons d\'eduit de la th\'eorie d'Arthur {\bf a minima} des limitations s\'ev\`eres sur les
caract\`eres infinit\'esimaux des repr\'esentations pouvant appara\^{\i}tre dans $L^2 (\Gamma \backslash G)$ lorsque $G^{{\rm nc}} = SU(n,1)$ ou $O(n,1)^0$. 
Commen\c{c}ons par rappeler cette th\'eorie d'Arthur {\bf a minima} et ce que celle-ci implique dans le cas des groupes de rang r\'eel $1$.

\paragraph{Param\`etres d'Arthur.} Dans ce paragraphe $G$ est un groupe r\'eductif r\'eel et ${}^L G$ le groupe dual de Langlands \cite{Langlands, Borel8}. Le {\it groupe
de Weil de ${\Bbb R}$}, not\'e $W_{{\Bbb R}}$, est l'extension de ${\Bbb C}^*$ par ${\Bbb Z} /2{\Bbb Z}$ (le groupe de Galois de ${\Bbb C}$ sur ${\Bbb R}$)~:
$$W_{{\Bbb R}} = {\Bbb C}^* \cup j {\Bbb C}^* ,$$
o\`u $j^2 = -1$ et $jcj^{-1} = \overline{c}$. Un {\it param\`etre d'Arthur pour $G$} est un homomorphisme 
\begin{eqnarray} \label{Aparam}
\psi : W_{{\Bbb R}} \times SL(2,{\Bbb C}) \rightarrow {}^L G
\end{eqnarray}
tel que 
\begin{enumerate}
\item Le diagramme 
$$\begin{array}{rcccl}
W_{{\Bbb R}} & & \rightarrow & & {}^L G \\
& \searrow &  & \swarrow & \\
& & {\rm Gal}({\Bbb C} /{\Bbb R}) && 
\end{array}$$
est commutatif.
\item Le morphisme $\psi_{| SL(2,{\Bbb C})}$ est holomorphe ($\equiv$ alg\'ebrique).
\item L'image de $\psi_{| W_{{\Bbb R}}}$ est d'adh\'erence compacte. 
\end{enumerate}
On d\'eduit de $\psi$ un ``param\`etre de Langlands''  (cf. \cite{Langlands, Borel8})
\begin{eqnarray} \label{AL}
\begin{array}{clcl}
\varphi_{\psi} : & W_{{\Bbb R}} & \rightarrow & {}^L G \\
& w & \mapsto & \psi \left( w, \left( 
\begin{array}{cc}
|w|^{1/2} & 0 \\
0 & |w|^{-1/2} 
\end{array} \right) \right)  
\end{array}
\end{eqnarray}
o\`u, pour $w \in W_{{\Bbb R}}$, $|w|$ est la valeur absolue de l'image de $w$ dans ${\Bbb R}^*$ via la suite exacte naturelle
$$1 \rightarrow U = \{ z \in {\Bbb C}^* \; : \; |z| =1 \} \stackrel{u}{\rightarrow} {\Bbb R}^* \rightarrow 1$$
donn\'ee par $u : z \mapsto |z|_{{\Bbb C}} = z \overline{z}$ si $z \in {\Bbb C}^*$ et $u(j)=-1$. 

Il correspond aux param\`etres d'Arthur des {\it paquets d'Arthur}, ensembles finis de repr\'esentations conjecturalement reli\'e \`a la 
d\'ecomposition du spectre automorphe. Nous avons \'ecrit ``param\`etre de Langlands'' entre guillemets car le param\`etre $\varphi_{\psi}$ ne d\'efinit en g\'en\'eral qu'une repr\'esentation de la forme int\'erieure quasi-d\'eploy\'ee du groupe $G$ (et non pas une repr\'esentation du groupe $G$ lui m\^eme). La d\'efinition des paquets d'Arthur
est donc un peu plus d\'elicates lorsque $G$ n'est pas quasi-d\'eploy\'e. Nous ne d\'etaillons pas cet aspect des choses, pr\'ef\'erant d'abord remarquer que 
puisque sur ${\Bbb C}$ les alg\`ebres de Lie de $G$ et de n'importe laquelle de ses formes
int\'erieures sont isomorphes, on peut en tout cas parler du {\it caract\`ere infinit\'esimal} d\'efinit par $\varphi_{\psi}$. \'Enon\c{c}ons alors la formulation tr\`es faible suivante des 
conjectures d'Arthur.

\begin{conj} \label{ArthurC}
Si une repr\'esentation irr\'eductible $\pi$ de $G$ appara\^{\i}t (faiblement) dans $L^2 (\Gamma \backslash G)$ pour un sous-groupe de congruence, son caract\`ere infinit\'esimal
est associ\'e \`a un param\`etre d'Arthur $\varphi_{\psi}$.
\end{conj}

Dans \cite{Livre} nous d\'eduisons de la Conjecture \ref{ArthurC} que lorsque $G^{{\rm nc}} = SU(n,1)$ ou $O(n,1)^0$, le groupe $G$ doit v\'erifier l'hypoth\`ese~:
\begin{eqnarray*}
({\rm Hyp}) & : & \left\{
\begin{array}{l}
\mbox{Si } \pi = \overline{X} (\lambda , s) \in \widehat{G}_{{\rm Aut}} \mbox{ est non temp\'er\'ee, alors} \\
s \in \frac12 {\Bbb Z} \; (0 < s \leq \rho ) .
\end{array}
\right.
\end{eqnarray*}
Ici, et comme au \S 1,  $\overline{X} (\lambda , s)$ d\'esigne l'unique quotient de Langlands de la repr\'esentation induite $X (\lambda , s)$ avec $\lambda$ (pseudo-)caract\`ere
r\'egulier associ\'e \`a une repr\'esentation du sous-groupe de Levi d'un parabolique minimal $P=MAN$ de $G$ et $s \in {\Bbb C}$ correspondant \`a un \'el\'ement de 
$\widehat{A}$ via l'identification de $1 \in {\Bbb C}$ avec l'unique racine simple de $A$ dans $G$. Enfin et de mani\`ere usuelle, $\rho = n$ si $G^{{\rm nc}} = SU(n,1)$
et $\rho = (n-1)/2$ si $G^{{\rm nc}} = O(n,1)^0$.

\bigskip

Revenons maintenant sur la th\'eorie d'Arthur dans le cas des repr\'esentations cohomologiques telle que d\'ecrite dans \cite{Arthur}. 
Commen\c{c}ons par d\'ecrire le foncteur d'induction cohomologique en termes de param\`etres de Langlands.
Consid\'erons donc $G$ un groupe semi-simple r\'eel et $L$ le sous-groupe de Levi associ\'e \`a une sous-alg\`ebre
parabolique $\theta$-stable $\mathfrak{q}$. Nous lui associons un plongement
$$\xi_L : {}^L L \rightarrow {}^L G .$$
La restriction de $\xi_L$ au groupe dual de $L$ est l'inclusion naturelle de celui-ci dans le groupe dual de $G$, alors que
$$\xi_L (z) = \left( \frac{z}{\overline{z}} \right)^{\rho (\mathfrak{u} )},$$
et
$$\xi_L (j) = n_G n_L^{-1} \rtimes j ,$$
o\`u $n_G$ (resp. $n_L$) est un \'el\'ement fix\'e du groupe d\'eriv\'ee du groupe dual de $G$ (reps. $L$) qui envoie les racines positives sur les racines
n\'egatives. Il d\'ecoule alors de \cite[Proposition 1.3.5]{Shelstad} que $\xi_L$ est un $L$-plongement. Il n'est de plus pas difficile de v\'erifier que si 
$\varphi : W_{{\Bbb R}} \rightarrow {}^L L$ est le param\`etre de Langlands de la repr\'esentation $Y = \overline{X}^L (\gamma_L )$ du point 4.
de la remarque suivant le Th\'eor\`eme \ref{rep isol} et v\'erifiant (\ref{conda}) et (\ref{condb}), alors $\xi_L \circ \varphi$ est le 
param\`etre de Langlands de la repr\'esentation ${\cal R}Y$. 

Le lien avec les param\`etres d'Arthur est le suivant. Consid\'erons le morphisme
$$\psi : W_{{\Bbb R}} \times SL(2,{\Bbb C}) \rightarrow {}^L L$$
dont la restriction \`a $W_{{\Bbb R}}$ est triviale et qui envoie l'\'el\'ement 
$$\left( 
\begin{array}{cc}
1 & 1 \\
0 & 1 
\end{array} \right)$$
sur l'\'el\'ement unipotent principal dans le groupe dual de $G$. C'est un param\`etre d'Arthur dont le param\`etre de Langlands induit $\varphi_{\psi}$ 
correspond \`a la repr\'esentation triviale de $L$. En composant le param\`etre $\psi$ par $\xi_L$, on obtient un param\`etre d'Arthur du groupe $G$, il lui 
correspond un paquet d'Arthur contenant la repr\'esentation cohomologique correspondant \`a la sous-alg\`ebre parabolique $\mathfrak{q}$.

Les conjectures d'Arthur reposent sur la construction de tels paquets, familles
finies de repr\'esentations de $G$ associ\'ees aux param\`etres d'Arthur, qui devraient partitionner le dual automorphe de $G$. 
Arthur lui-m\^eme ne dit rien de la construction de ces paquets, sauf dans le cas des repr\'esentations cohomologiques, via la $L$-application $\xi_L$ construite ci-dessus. En particulier, sa formulation ne pr\'ecise pas quels param\`etres de Langlands devraient 
appara\^{\i}tre dans la paquet associ\'e \`a un param\`etre d'Arthur donn\'e, quand $G$ n'est pas quasi-d\'eploy\'e sur ${\Bbb R}$. 
Ces param\`etres sont a priori construits par Adams, Barbasch et Vogan \cite{AdamsBarbaschVogan}.
Mais la construction qui repose sur la g\'eom\'etrie alg\'ebrique, les singularit\'es d'orbites et les faisceaux pervers est tr\`es difficile. 
Notre interpr\'etation, plus modeste, des conjectures d'Arthur dans ce contexte implique la conjecture suivante.
Toute erreur est \'evidemment notre. 

\begin{conj} \label{conj gal}
Soit $G$ un groupe semi-simple alg\'ebrique sur ${\Bbb Q}$ et $L$ le sous-groupe de Levi associ\'e \`a une sous-alg\`ebre
parabolique $\theta$-stable $\mathfrak{q}$. Fixons $\tau$ une repr\'esentation unitaire de $L$ et $\pi$ la repr\'esentation de $G$ d\'eduite de $\tau$ par 
la fonctorialit\'e induite par la $L$-application $\xi_L$. La repr\'esentation $\pi$ appartient au dual automorphe $\widehat{G}_{{\rm Aut}}$ de $G$ si et seulement si
(quitte \`a conjuguer $L$ dans $G$) l'inclusion $L\subset G$ est d\'efinie sur ${\Bbb Q}$ et la repr\'esentation $\tau \in \widehat{L}_{{\rm Aut}}$.
\end{conj}

L'induction cohomologique devrait donc \^etre fonctorielle entre les duaux automorphes. 
Le lien avec la Conjecture \ref{isolaut} est le suivant. D'apr\`es la Remarque 5 du \S 1, 
il d\'ecoule de la d\'emonstration originale par Vogan \cite{Vogan} du Th\'eor\`eme \ref{rep isol} que pour 
d\'emontrer la Conjecture \ref{isolaut} il suffit de v\'erifier que si $\pi$ est comme dans le Th\'eor\`eme \ref{rep isol}, que $L$ a un centre compact et que 
$\tau_i$ est une suite de repr\'esentations unitaires de $L$ qui converge vers $1_L$, alors pour $i$ suffisamment grand ${\cal R} \tau_i$ n'appartient
pas au dual automorphe de $G$. Puisque $L$ a la propri\'et\'e $\tau$, cf. Clozel \cite{Clozel}, la Conjecture \ref{conj gal} implique finalement la 
Conjecture \ref{isolaut}.

\medskip

Plus modestement nous montrerons au \S 5 comment, dans de nombreux cas, la Conjecture \ref{isolaut} se ram\`ene 
\`a v\'erifier l'hypoth\`ese $({\rm Hyp})$, ou une approximation de celle-ci, pour chacun des groupes
$O(n,1)^0$ $(n\geq 2)$ ou $SU(n,1)$ $(n\geq 1)$. Nous en d\'eduirons de nouveaux cas de la Conjecture \ref{isolaut}.

Avant cela et dans la section suivante, nous donnons d'autres raisons de croire \`a la Conjecture \ref{conj gal}. 

\section{Applications cohomologiques}

Consid\'erons toujours $G$ un groupe semi-simple alg\'ebrique sur ${\Bbb Q}$. 
Les applications cohomologiques que nous avons en vus sont des r\'eponses partielles \`a la question classique suivante~:
\begin{Q} \label{Q1}
Soit $\pi \in \widehat{G}$ une repr\'esentation cohomologique. Existe-t-il un sous-groupe de congruence
$\Gamma \subset G({\Bbb Q})$ tel que la repr\'esentation $\pi$ intervienne discr\`etement dans $L^2 (\Gamma \backslash G )$ ?
\end{Q}
ou \`a la question plus faible~:
\begin{Q} \label{Q2}
Soit  $\pi \in \widehat{G}$ une repr\'esentation cohomologique. Existe-t-il un r\'eseau $\Gamma$ dans $G$ tel que la repr\'esentation $\pi$ intervienne discr\`etement dans
$L^2 (\Gamma \backslash G )$ ?
\end{Q}

Supposons $G$ d\'efinit sur ${\Bbb Q}$. Notre strat\'egie, qui est d\'ej\`a celle de Burger et Sarnak \cite{BurgerSarnak}, 
pour montrer qu'une repr\'esentation cohomologique $\pi$ intervient discr\`etement dans
$L^2 (\Gamma \backslash G )$ pour un certain sous-groupe de congruence de $\Gamma \subset G({\Bbb Q})$ consiste 
1) \`a d\'emontrer que la repr\'esentation $\pi$ appartient au dual automorphe $\widehat{G}_{{\rm Aut}}$, et 2) \`a v\'erifier qu'elle y est isol\'ee.

\medskip

Commen\c{c}ons par remarquer que la Conjecture \ref{conj gal} implique que pour tout sous-groupe de Levi
$L \subset G$ d\'efini sur ${\Bbb Q}$ et associ\'e \`a une sous-alg\`ebre
parabolique $\theta$-stable $\mathfrak{q} \subset \mathfrak{g}$, la
repr\'esentation cohomologique ${\cal R} 1_L  \in \widehat{G}_{{\rm
Aut}}$. En particulier et en utilisant la Conjecture \ref{isolaut}, 
la Conjecture \ref{conj gal} implique~:
\begin{conj}
Soit $G$ un groupe semi-simple qui admet une s\'erie discr\`ete. Alors la r\'eponse \`a la Question \ref{Q2} est toujours positive.
\end{conj}
La Question \ref{Q1} est plus d\'elicate. Un bon exemple est celui des groupes unitaires tordus \'etudi\'es par Clozel \cite{Clozel2}. Afin de motiver
la Conjecture \ref{conj gal}, montrons comment d\'eduire de celle-ci les
principaux r\'esultats de \cite{Clozel2}.

Fixons trois entiers $n$, $p$ et $q$ tels que $n=p+q$, $1\leq p\leq q$. Si
$b$ est un diviseur de $n$, rappelons que Clozel pose
\begin{eqnarray}
x=\left[ \frac{p}{b} \right]
\end{eqnarray}
o\'u $[t]$ est la partie enti\`ere de $t\geq0$;
\begin{eqnarray}
n=ab;
\end{eqnarray}
\begin{eqnarray}
N=N(b) = bx^2 + (b-2p) x + (a-1) p .
\end{eqnarray}
Enfin, si $N=N(b)$ est un tel entier, on note $[pq-N,pq+N]_2$ l'ensemble
des entiers $\{ pq-N , pq-N+2 , \ldots , pq+N \}$ de m\^eme parit\'e.

\begin{prop} \label{clozel}
Soient $F$ un corps de nombres totalement r\'eel, $F_c$ une extension de $F$ quadratique et totalement imaginaire (corps CM).
Soit $(D,*)$ une alg\`ebre \`a division de degr\'e $n^2$ sur $F_c$, munie d'une involution de seconde
esp\`ece $*$. Soient $U$ le groupe unitaire sur $F$ associ\'ee \`a $(D,*)$ et $G={\rm Res}_{F/{\Bbb Q}} (U/F)$ - groupe ${\Bbb Q}$-alg\'ebrique.
Supposons le groupe unitaire $U$ compact en toutes les places infinies de $F$ sauf une.
On a alors $G^{{\rm nc}} = SU(p,q)$, pour certains entiers $p$, $q$ tels que $p+q=n$, et la Conjecture \ref{conj gal} implique~:
\begin{enumerate}
\item Une repr\'esentation non triviale $\pi \in \widehat{G}_{{\rm Aut}}$ ne peut intervenir dans la cohomologie qu'en les degr\'es appartenant \`a
la r\'eunion, sur l'ensemble des diviseurs $b\neq1$ de $n$, des intervalles
$$[pq-N(b) , pq+N(b) ]_2 .$$
\item Pour chaque degr\'e du point 1. il existe un sous-groupe de
congruence $\Gamma \subset G({\Bbb Q} )$ et une repr\'esentation
$\pi \subset L^2 (\Gamma \backslash G )$ cohomologique en ce degr\'e.
\end{enumerate}
\end{prop}
{\it D\'emonstration.} Nous adoptons ici un style volontairement l\'eger, notre but est en effet de motiver la Conjecture \ref{conj gal}. 
Qu'il nous suffise donc de remarquer que les arguments de cohomologie galoisienne
d\'evelopp\'es dans \cite{Clozel} montrent que la seule mani\`ere de plonger
un groupe unitaire $H$ comme Levi dans $G$ est de consid\'erer un diviseur $b$ de $n$,
une extension $L_0$ de $F$, totalement r\'eelle, de degr\'e $b$ et une
alg\`ebre \`a division $A$ sur $L=L_0 \otimes F_c$ de degr\'e $a=n/b$ et
incluse dans $D \otimes L_0$ telle que, munie de l'involution de seconde
esp\`ece (relativement \`a $L/L_0$) induite par $*$, le groupe $U(A,*)$ -
un groupe sur $L_0$ - se plonge dans $U(D,*)$. Le groupe $H$ s'obtient
alors par restriction des scalaires \`a partir du groupe $U(A,*)$. On v\'erifie
ais\'ement que $H^{{\rm nc}}$ est isomorphe \`a un groupe du type 
\begin{eqnarray} \label{L}
SU(x_1 , y_1 ) \times \ldots \times SU(x_b , y_b ),
\end{eqnarray}
avec $x_i + y_i = a$ ($i=1, \ldots ,b$), $x_1 + \ldots + x_b = p$ et $y_1 +
\ldots + y_b = q$.

La Conjecture \ref{conj gal} implique que les seules repr\'esentations
cohomologiques qui apparaissent dans le dual automorphe
$\widehat{G}_{{\rm Aut}}$ sont obtenues par induction cohomologique de la
repr\'esentation triviale d'un ${\Bbb Q}$-sous-groupe $L\subset G$. Un tel
sous-groupe est n\'ecessairement de la forme (\ref{L}) et la
repr\'esentation cohomologique correspondante n'intervient que dans les
degr\'es appartenant \`a l'intervalle
$$[pq-\sum_{i=1}^b x_i y_i  , pq +\sum_{i=1}^b x_i y_i]_2 .$$
Le premier point de la Proposition d\'ecoule donc du petit lemme
combinatoire suivant.              
\begin{lem} \label{lemC}
Supposons que $p+q=ab$ et fixons $b$ entiers $x_1 , \ldots , x_b$ (resp. $y_1 , \ldots , y_b$) tels que
\begin{itemize}
\item $x_1 + \ldots + x_b=p$, $y_1 + \ldots + y_b=q$, et
\item pour tout $i=1, \ldots , b$, $x_i+y_i = a$.
\end{itemize}
Alors,
$$\sum_{i=1}^b x_i y_i \leq (b-y)x(a-x) + y(x+1)(a-x-1),$$
o\`u $p=bx+y$ avec $0 \leq y <b$ et ces deux nombres ont la m\^eme parit\'e. Autrement dit
$$\sum_{i=1}^b x_i y_i \leq N(b), $$
et ces deux nombres ont la m\^eme parit\'e. 
\end{lem}
{\it D\'emonstration.} Il s'agit en fait d'\'evaluer la somme 
\begin{eqnarray} \label{s}
\sum_{i=1}^b x_i (a-x_i) 
\end{eqnarray}
lorsque les $x_i$ sont des entiers astreints \`a appartenir \`a
l'intervalle $[0,a]$ et \`a v\'erifier l'identit\'e $x_1 + \ldots + x_b =
p$. Puisque $(x-1)(a-x+1) = x(a-x) -a-1-2x$, il est imm\'ediat, que sous ces
conditions changer les $x_i$ ne change pas la parit\'e de l'expression
(\ref{s}). Il nous reste \`a majorer (\ref{s})~: on montre par r\'ecurrence
que l'expression (\ref{s}) est maximale lorsque $x_1 = \ldots = x_{b-y} =
x$ et $x_{b-y+1} = \ldots = x_b = x+1$. Il suffit de traiter le cas $b=2$
qui est bien s\^ur imm\'ediat. Enfin il n'est pas difficile de v\'erifer
que $N(b) =(b-y)x(a-x) + y(x+1)(a-x-1)$. Le Lemme \ref{lemC} est d\'emontr\'e.

\medskip

Il nous reste \`a expliquer pourquoi en admettant la Conjecture \ref{conj
gal}, les degr\'es du point 1. interviennent effectivement dans la
cohomologie. Cela d\'ecoule de ce que l'on peut effectivement construire
comme indiqu\'e ci-dessus (cf. l\`a encore, les techniques de \cite{Clozel}) un
${\Bbb Q}$-sous-groupe $L\subset G$ tel que 
$$L^{{\rm nc}} = SU(x,a-x)^{b-y} \times SU(x+1 , a-x-1)^y.$$
(Notations du Lemme \ref{lemC}.)
La Conjecture \ref{conj gal} implique alors que la repr\'esentation $\pi$ obtenue
par induction cohomologique de $L$ \`a $G$ de la repr\'esentation triviale
appartient au dual automorphe de $G$ ($L$ est bien le sous-groupe de Levi
d'une sous-alg\`ebre parabolique $\theta$-stable de $\mathfrak{g}$). Mais
d'apr\`es la Conjecture \ref{isolaut} (qui d\'ecoule elle-m\^eme,
rappelons-le, de la Conjecture \ref{conj gal}), les repr\'esentations
cohomologiques sont isol\'ees dans $\widehat{G}_{{\rm Aut}}$. La
repr\'esentation $\pi$ doit donc intervenir discr\`etement dans
$\widehat{G}_{{\rm Aut}}$ et donc dans $L^2 (\Gamma \backslash G)$ pour un
certain sous-groupe de congruence $\Gamma \subset G({\Bbb Q})$. Un petit
calcul montre enfin que la repr\'esentation $\pi$ intervient dans la
cohomologie pour tous les degr\'es $\in [pq-N(b) , pq+N(b)]_2$. La
Proposition \ref{clozel} est d\'emontr\'ee.

\medskip

\noindent
{\bf Remarques.} {\bf 1.} Sous l'hypoth\`ese technique suppl\'ementaire (not\'ee (R) dans \cite{Clozel2})~:
\begin{center}
en chaque place $v$ de $F_c$, l'alg\`ebre $D_v = D\otimes_{F_c} F_{c,v}$ est soit
isomorphe \`a $M_n (F_{c,v} )$ soit une alg\`ebre \`a division,
\end{center}
la conclusion de la Proposition \ref{clozel} est un th\'eor\`eme de Clozel \cite[Theorem 3.3]{Clozel2}.

{\bf 2.} Nous pourrions traiter de mani\`ere analogue tous les groupes unitaires, ceux-ci sont plus g\'en\'eralement obtenus \`a partir d'une alg\`ebre 
centrale simple $B=M_r (D)$, avec $D$ comme dans l'\'enonc\'e de la Proposition \ref{clozel}.

\bigskip

Dans un m\^eme registre, remarquons qu'il n'est pas difficile de v\'erifier, cf. \cite{PlatonovRapinchuk},
que si $G$ est un groupe semi-simple alg\'ebrique sur ${\Bbb Q}$ provenant (par restriction des scalaires) d'un groupe de type
${}^3 D_4$ ou ${}^6 D_4$ sur un corps de nombre totalement r\'eel et tel que $G^{{\rm nc}} = O(7,1)^0$, alors $G$ ne contient aucun ${\Bbb Q}$-sous-groupe $L$ tel que 
$L^{{\rm nc}} = O(5,1)^0$. Et la Conjecture \ref{conj gal} implique la conjecture suivante.

\begin{conj} \label{nonTh}
Soit $G$ un groupe semi-simple alg\'ebrique sur ${\Bbb Q}$ provenant (par restriction des scalaires) d'un groupe de type
${}^3 D_4$ ou ${}^6 D_4$ sur un corps de nombre totalement r\'eel et tel que $G^{{\rm nc}} = O(7,1)^0$, alors pour tout sous-groupe 
de congruence $\Gamma \subset G({\Bbb Q})$, on a~:
$$b_1 (\Gamma ) =0 .$$
\end{conj}

Rappelons qu'une c\'el\`ebre conjecture de Thurston affirme que toute vari\'et\'e hyperbolique compacte admet un rev\^etement fini avec un premier
nombre de Betti non nul. Cette conjecture est v\'erifi\'ee pour toute vari\'et\'e hyperbolique arithm\'etique de dimension $\geq 4$ associ\'ee \`a un groupe $G$, $\neq {}^{3,6} D_4$, en restant
dans le monde des sous-groupes de congruence. Il est surprenant de noter que selon la Conjecture \ref{nonTh}, ce ne devrait plus \^etre le cas pour les groupes de type ${}^{3,6} D_4$.

\bigskip

\`A d\'efaut de savoir d\'emontrer la Conjecture \ref{conj gal}, on peut utiliser les fonctorialit\'es
remarquables d\'ecouvertes par Burger, Li et Sarnak dans \cite{BurgerLiSarnak} et \cite{BurgerSarnak}~:
\begin{enumerate}
\item Si $H \subset G$ sont deux ${\Bbb Q}$-groupes semi-simples et si $\rho \in \widehat{H}_{{\rm Aut}}$
alors toute repr\'esentation dans le support de l'induite ${\rm ind}_H^G \rho$ appartient au dual automorphe $\widehat{G}_{{\rm Aut}}$.
\item Si $H \subset G$ sont deux ${\Bbb Q}$-groupes semi-simples et si $\pi \in \widehat{G}_{{\rm Aut}}$
alors toute repr\'esentation dans le support de la restriction $\pi_{|H}$ appartient au dual automorphe $\widehat{H}_{{\rm Aut}}$.
\end{enumerate}

Commen\c{c}ons par exploiter la premi\`ere fonctorialit\'e en prenant pour $\rho$ la repr\'esentation triviale $1_H$ de $H$. Alors 
${\rm ind}_H^G \rho = {\rm ind}_H^G 1_H = L^2 (H \backslash G)$. Toute repr\'esentation de la s\'erie discr\`ete de $H \backslash G$ 
({\it i.e.} toute repr\'esentation intervenant discr\`etement dans la repr\'esentation
$L^2 (H \backslash G)$) appartient au dual automorphe $\widehat{G}_{{\rm Aut}}$. Comme le remarquent Burger, Li et Sarnak, une repr\'esentation cohomologique $\pi$ {\bf isol\'ee} de $G$
intervenant dans la s\'erie discr\`ete de $H \backslash G$ pour un certain sous-groupe $H$ comme ci-dessus intervient alors dans la cohomologie ($L^2$) de $S(\Gamma) $ pour un certain
sous-groupe de congruence $\Gamma \subset G({\Bbb Q})$. Nous montrons en fait la proposition suivante.

\begin{prop} \label{P1}
Soit $\pi$ une repr\'esentation cohomologique de $G$ appartenant \`a la s\'erie discr\`ete de $H \backslash G $.
Supposons de plus $\pi$ {\bf isol\'ee sous la condition $d=0$} dans $\widehat{G}_{{\rm Aut}}$.
Il existe alors un sous-groupe de congruence $\Gamma \subset G({\Bbb Q})$ tel que la repr\'esentation $\pi$
intervienne discr\`etement dans $L^2 (\Gamma \backslash G )$.
\end{prop}
{\it D\'emonstration.} Notons $R$ le degr\'e fortement primitif de la repr\'esentation cohomologique $\pi$.
L'espace ${\rm Hom}_K (\bigwedge^R \mathfrak{p} , \pi )$ est de dimension $1$, notons $\omega$ un g\'en\'erateur.
L'\'el\'ement $\omega \in {\rm Hom}_K (\bigwedge^R \mathfrak{p} , C^{\infty} (H \backslash G))$ et d\'efinit donc
une forme diff\'erentielle lisse sur $X_G$, $H$-invariante \`a gauche. La repr\'esentation $\pi$ \'etant cohomologique, 
cette forme diff\'erentielle est de plus harmonique. Elle s'identifie naturellement \`a une composante isotypique 
de la forme de Thom $\Omega$ (\`a une forme exacte lisse pr\`es) du fibr\'e normal de $X_H$ dans $X_G$, cf. \cite{TongWang}.  
Or il existe (cela d\'ecoule de la construction de la forme de Thom) une forme diff\'erentielle $\alpha$ sur $X_G$ lisse sauf peut-\^etre le long de 
$X_H$ et telle que $\Omega = d(\alpha )$. Puisque le fibr\'e normal comme la forme $\Omega$ sont $H$-invariants \`a gauche, on peut
de plus supposer que la forme $\alpha$ est $H$-invariante \`a gauche et la repr\'esenter par un \'el\'ement
$$\alpha \in  {\rm Hom}_K (\bigwedge {}^{R-1} \mathfrak{p} , C^{\infty} (H \backslash G)).$$

Fixons maintenant $d$ une distance $G$-invariante \`a gauche sur $G$. Notons 
$$\varphi (x) = \left\{ 
\begin{array}{cl}
e^{\frac{1}{x^2 -1}} & \mbox{ si } |x|<1 \\
0 & \mbox{ si } |x| \geq 1
\end{array} \right.$$
($x \in {\Bbb R}$) et $C= \left( \int \varphi \right)^{-1}$. Consid\'erons alors pour $s>0$, la forme diff\'erentielle lisse \`a support compact~:
$$\Omega_s = Csd\left[(\varphi (sx) * 1_{[0,s]}) (d(\cdot , H)) \cdot \alpha \right] \in {\rm Hom}_K (\bigwedge {}^R \mathfrak{p} , C_0^{\infty} (H \backslash G)).$$
Les formes $\Omega_s$ convergent uniform\'ement sur les compacts vers la forme $\Omega$ lorsque $s$ tend vers $+\infty$. 

\'Etant donn\'e un sous-groupe de congruence $\Gamma \subset G({\Bbb Q})$, on peut toujours former la s\'erie
$$\Omega^{\Gamma}_s = \sum_{\gamma \in (\Gamma \cap H) \backslash \Gamma} \gamma^* \Omega_s ,$$
qui d\'efinit une forme diff\'erentielle sur $S(\Gamma)$. Fixons $\Gamma$ et notons $\Lambda : = \Gamma \cap H$.
Il est bien connu, cf. \cite{Livre}, qu'il existe une suite $\Gamma_m $ de sous-groupes de congruence $\subset \Gamma$ telle que
$$\Lambda = \cap_m \Gamma_m .$$
Les formes $\Omega^{\Gamma_m}_s$ convergent donc sur tout compact de $X_G$ vers la forme $\Omega_s$. La suite de formes diff\'erentielles {\bf ferm\'ees}
$\Omega^{\Gamma_m}_{m}$ converge donc vers $\Omega$. Mais la repr\'esentation $\pi$ (qui intervient dans $\Omega$) est isol\'ee sous la condition $d=0$ 
dans le dual automorphe de $G$, elle doit donc n\'ecessairement intervenir dans $\Omega^{\Gamma_m}_{m}$ pour $m$
suffisamment grand et donc dans le spectre discret de $\Gamma_m \backslash G$ pour $m$ suffisamment grand. Ce qui conclut la d\'emonstration 
de la Proposition \ref{P1}.

\bigskip

De mani\`ere analogue on montre la proposition suivante dans \cite{IRMN}.

\begin{prop} \label{P2}
Soit $\pi$ (resp. $\sigma$) une repr\'esentation cohomologique de $G$ (resp. $H$) de degr\'e fortement primitif $R$. Notons
$\mu$ et $\delta$ leurs plus bas $K$-types respectifs. Supposons que 
\begin{enumerate}
\item $\sigma$ (resp. $\delta$) intervienne discr\`etement dans la restriction de $\pi$ (resp. $\mu$) \`a $H$ (resp. $K\cap H$),
\item $\pi \in \widehat{G}_{{\rm Aut}}$, et
\item $\sigma$ soit {\bf isol\'ee sous la condition $d=0$} dans $\widehat{H}_{{\rm Aut}}$. 
\end{enumerate}
Alors il existe un sous-groupe de congruence $\Lambda \subset H({\Bbb Q})$ tel que la repr\'esentation $\sigma$ intervienne 
discr\`etement dans $L^2 (\Lambda \backslash H )$.
\end{prop}

\bigskip

Nous d\'eduisons maintenant de ces Propositions et de la Proposition \ref{d} des r\'eponses partielles aux Questions \ref{Q1} et \ref{Q2}.
De la Proposition \ref{P1} nous d\'eduisons d'abord le th\'eor\`eme g\'en\'eral suivant.

\begin{thm} \label{T}
Soient $G$ et $H$ deux groupes semi-simples tels que $H$ co\"{\i}ncide avec 
l'ensemble des points fixes d'une involution de $G$ qui commute \`a l'involution de Cartan. Supposons~:
\begin{enumerate}
\item qu'il existe une sous-alg\`ebre de Cartan {\bf compacte} $\mathfrak{t}$ du suppl\'ementaire de $\mathfrak{h}$ dans $\mathfrak{g}$ \footnote{Ce qui est \'equivalent au fait que $H\backslash G$ poss\`ede une s\'erie discr\`ete.}, et
\item qu'il n'existe aucune sous-alg\`ebre parabolique $\theta$-stable $\mathfrak{q}' = \mathfrak{u}' + \mathfrak{l}'$ telle que le groupe $L'$
ait un facteur simple localement isomorphe \`a $SO(n,1)$ ($n\geq 2$) ou \`a $SU(n,1)$ ($n\geq 1$) et que le centralisateur $L$ de $\mathfrak{t}$ dans $G$
s'obtienne en rempla\c{c}ant ce facteur par $SO(n-2,1) \times SO(2)$ ou $S(U(n-1, 1) \times U(1))$ respectivement.
\end{enumerate} 
Il existe alors un r\'eseau cocompact $\Gamma \subset G$ tel que 
$$H^{d_G - d_H} (S(\Gamma )) \neq 0 ,$$
o\`u $d_G$ (resp. $d_H$) d\'esigne la dimension de l'espace sym\'etrique associ\'e \`a $G$ (resp. $H$). 
Il existe plus pr\'ecis\'ement un r\'eseau cocompact $\Gamma \subset G$ tel que la partie 
$H^{d_G - d_H} (S(\Gamma ))_{2(\rho-\rho_c)}$, de la cohomologie, qui correspond au $K$-type $2(\rho -\rho_c)$, soit {\bf non nulle}.
\end{thm}
{\it D\'emonstration.} Quitte \`a ajouter des facteurs compacts \`a $G$ et $H$, nous pouvons les supposer alg\'ebriques sur ${\Bbb Q}$, avec $G({\Bbb Q})$ cocompact dans $G({\Bbb A})$, 
et supposer le plongement $H\subset G$ rationnel. Nous pouvons enfin supposer le groupe $G^{{\rm nc}}$ simple.

Notons $\pi$ la repr\'esentation cohomologique de $G$ de plus bas $K$-type $2(\rho - \rho_c )$, remarquons que son 
degr\'e fortement primitif est \'egal \`a $d_G - d_H$. La repr\'esentation $\pi$ intervient dans la s\'erie discr\`ete de $H\backslash G$ (cf. par exemple \cite{Lefschetz})~:
elle est associ\'ee \`a une sous-alg\`ebre parabolique $\theta$-stable $\mathfrak{q} \subset \mathfrak{g}$ dont la sous-alg\`ebre de Levi $\mathfrak{l}$ est \'egale au centralisateur de
$\mathfrak{t}$ dans $\mathfrak{g}$.   

D'apr\`es la Proposition \ref{d}, l'hypoth\`ese 2. du Th\'eor\`eme \ref{T} implique que la repr\'esentation $\pi$ est {\bf isol\'ee sous la condition $d=0$} dans $\widehat{G}$. Le 
Th\'eor\`eme \ref{T} d\'ecoule donc  de la Proposition \ref{P1}.

\bigskip
 
Le Th\'eor\`eme \ref{T} s'applique notamment aux couples de groupes suivants~:
$(H,G)=(S(U(p,q-r) \times U(r)) , SU(p,q))$, $(S(O(p,q-r) \times O(r))^0 , O(p,q)^0 )$ ou $(Sp (p,q-r) \times Sp(r) , Sp(p,q))$. Les hypoth\`eses du Th\'eor\`eme \ref{T} sont v\'erifi\'ees dans les deux premiers cas dans \cite{Lefschetz},
le dernier cas est similaire. Nous obtenons le corollaire suivant.

\begin{cor} \label{ttt}
\begin{enumerate}
\item Soit $G=SU(p,q)$ avec $p\geq 2$. Alors la r\'eponse \`a la Question \ref{Q2} est positive pour chacune des repr\'esentations
$A((r^p ) , ((q-r)^p))$, $0 \leq r \leq q/2$.
\item Soit $G=O(p,q)^0$ avec $p\geq 2$. Alors la r\'eponse \`a la Question \ref{Q2} est positive pour chacune des repr\'esentations
$A((r^p ))$, $0 \leq r \leq q/2$.
\item Soit $G=Sp(p,q)$ avec $p\geq 2$. Alors la r\'eponse \`a la Question \ref{Q2} est positive pour chacune des repr\'esentations
$A((0) , ((q-2r)^p))_0$, $0 \leq r \leq q/2$.
\end{enumerate}
\end{cor}

\medskip

\paragraph{Rel\`evement th\'eta.} La m\'ethode la plus puissante pour s'attaquer aux Questions \ref{Q1} et \ref{Q2} semble \^etre la th\'eorie du relev\'e th\'eta. Rappelons bri\`evement comment
celle-ci fonctionne. Soit $k$ un corps de nombre totalement r\'eel et soit $(D,\i)$ une $k$-alg\`ebre \`a involution de l'un des trois types suivants~:
\begin{eqnarray} \label{D}
D = \left\{
\begin{array}{ll}
k, & \mbox{cas 1}, \\
\mbox{une extension quadratique } F/k, & \mbox{cas 2}, \\
\mbox{une alg\`ebre de quaternion de centre } k, & \mbox{cas 3}, 
\end{array} \right. 
\end{eqnarray} 
et
\begin{eqnarray} \label{i}
\i = \left\{
\begin{array}{ll}
id, & \mbox{cas 1}, \\
\mbox{l'involution de Galois de } F/k, & \mbox{cas 2}, \\
\mbox{l'involution standard}, & \mbox{cas 3}. 
\end{array} \right. 
\end{eqnarray} 
Soient $V$ et $V'$ deux espaces vectoriels de dimension finie sur $D$ \'equip\'es de deux formes sesquilin\'eaires non-d\'eg\'en\'er\'ees $(.,.)$ et $(.,.)'$, l'une $\i$-hermitienne, l'autre $\i$-anti-hermitienne.
Le $k$-espace vectoriel $W=V\otimes_D V'$ est alors naturellement muni d'une forme symplectique
$$\langle .,. \rangle = {\rm tr}_{D/k} ((.,.) \otimes (.,.)' {}^{\i} )$$
o\`u ${\rm tr}_{D/k}$ d\'esigne la trace usuelle de $D$ sur $k$. Notons $Sp(W)$, $G$ et $G'$ les groupes d'isom\'etries respectifs de $\langle .,. \rangle$, $(.,.)$ et $(.,.)'$. Alors 
$(G,G')$ es une paire r\'eductive duale irr\'eductible de type I dans $Sp(W)$. Supposons le groupe $Sp(W)$ compact \`a toutes les places \`a l'infini de $k$ sauf une et notons 
$Sp$ le groupe r\'eel (non compact) en cette place. Nous noterons \'egalement $G$ et $G'$ les sous-groupes r\'eels de $Sp$ correspondants.
Notons enfin $\widetilde{Sp}$ le rev\^etement m\'etaplectique \`a deux feuillets du groupe symplectique $Sp$ et 
$\widetilde{G}$, $\widetilde{G}'$ les images inverses respectives de $G$ et $G'$ dans $\widetilde{Sp}$. Comme dans l'introduction, et apr\`es restriction des scalaires de 
$k$ \`a ${\Bbb Q}$, nous pouvons parler du dual automorphe de $\widetilde{Sp}$. Traduits en ces termes, Howe montre, entre autres choses, dans
\cite{Howe}, que la repr\'esentation de Weil (ou repr\'esentation de l'oscillateur harmonique) $\omega$ du groupe $\widetilde{Sp}$ appartient au dual automorphe de $\widetilde{Sp}$. 
Consid\'erons maintenant deux repr\'esentations respectives $\pi$ et $\pi'$ de $G$ et $G'$. Les deux repr\'esentations $\pi$ et $\pi'$ sont dites duales pour la correspondance 
$\theta$ locale, not\'e $\pi' = \theta (\pi)$, si la repr\'esentation $\pi \otimes \pi '$ de $\widetilde{G} \cdot \widetilde{G}'$ est \'equivalente \`a une sous-repr\'esentation irr\'eductible de la restriction de $\omega$
\`a $\widetilde{G} \cdot \widetilde{G}'$. Dans \cite{Li}, Li montre que si $\pi$ est une repr\'esentation de la s\'erie discr\`ete de $G$ ``suffisamment r\'eguli\`ere'' (cf. \cite{Li}) et si
${\rm dim}_D V \geq {\rm dim}_D V'$, alors $\pi'$ admet un relev\'e th\'eta non nul \`a $\widetilde{G}$. La repr\'esentation correspondante $\pi = \theta (\pi ')$ est 
une repr\'esentation cohomologique explicitement d\'ecrite dans \cite{Li}.

Dans \cite{Li2} Li v\'erifie que dans un grand nombre de cas la correspondance est globale. Ce qui lui permet de d\'emontrer le th\'eor\`eme suivant (le groupe $G$ est toujours comme ci-dessus).

\begin{thm} \label{Li gen}
\begin{enumerate}
\item Soit $G=SU(p,q)$. Alors la r\'eponse \`a la Question \ref{Q2} est positive pour chacune des repr\'esentations
$A(\lambda , \mu)$, avec $\mu/ \lambda$ rectangle de p\'erim\`etre $>p+q$.
\item Soit $G=O(p,q)^0$ avec $p+q$ pair. Alors la r\'eponse \`a la Question \ref{Q1} est positive pour chacune des repr\'esentations
$A(\lambda)$, avec $\hat{\lambda}/ \lambda$ rectangle de p\'erim\`etre $>p+q +2$.
\item Soit $G=O(p,q)^0$ avec $p+q$ impair. Alors la r\'eponse \`a la Question \ref{Q1} est positive pour chacune des repr\'esentations
$A((r^p))$, avec $r<\frac{1}{4} (p+q-2)$.
\item Soit $G=Sp(p,q)$. Alors la r\'eponse \`a la Question \ref{Q1} est positive pour chacune des repr\'esentations
$A(\lambda , \mu)_0$, avec $\mu/ \lambda$ rectangle de p\'erim\`etre $\geq p+q$.
\end{enumerate}
\end{thm}

Dans le cas des espaces hermitiens et concernant la cohomologie {\bf holomorphe} Anderson \cite{Anderson}
r\'epond positivement \`a la Question \ref{Q2} en utilisant l\`a encore la technique du rel\`evement th\'eta.
Remarquons que dans le Th\'eor\`eme \ref{Li gen} le cas du groupe unitaire est sp\'ecial puisque seule la Question \ref{Q2} est r\'esolue pour les repr\'esentations
cohomologiques consid\'er\'ees. Ceci ne doit pas nous \'etonner, la r\'eponse \`a la Question \ref{Q1} est en effet n\'egative, cf. Proposition \ref{clozel} ou mieux \cite{Clozel2}.

\medskip

Lorsque la repr\'esentation cohomologique est isol\'ee dans le dual automorphe, on peut directement d\'eduire de la correspondance locale
la correspondance globale. Il d\'ecoule en effet de la deuxi\`eme fonctorialit\'e de Burger, Li et Sarnak que si la repr\'esentation $\pi \otimes \pi'$
de $\widetilde{G} \cdot \widetilde{G}'$ est \'equivalente \`a une sous-repr\'esentation irr\'eductible de la restriction de la repr\'esentation $\omega \in \widehat{(\widetilde{Sp})}_{{\rm Aut}}$, alors
$\pi \in \widehat{G}_{\rm Aut}$ (et $\pi' \in \widehat{G'}_{\rm Aut}$). Si $\pi$ est isol\'ee dans $\widehat{G}_{\rm Aut}$, elle doit donc intervenir discr\`etement dans un $L^2 (\Gamma \backslash G)$, 
pour $\Gamma$ sous-groupe de congruence de $G$. Le th\'eor\`eme et la conjecture qui suivent d\'ecoulent donc des travaux \cite{Li} de Li et respectivement du Th\'eor\`eme \ref{rep isol} et de la Conjecture \ref{isolaut}.

\begin{thm} \label{relth}
Soit $G$ comme au-dessus et soit $\pi = A_{\mathfrak{q}}$ une repr\'esentation cohomologique de $G$ telle que
$$\mathfrak{l}_0 \cong \mathfrak{l}_0' \oplus \mathfrak{g}_1$$
o\`u
\begin{enumerate}
\item $ \mathfrak{l}_0' $ est une alg\`ebre de Lie compacte;
\item $\mathfrak{g}_1$ est du ``m\^eme type'' que $\mathfrak{g}_0$, c'est-\`a-dire que $\mathfrak{g}_1$ est isomorphe \`a l'alg\`ebre de Lie du groupe
des isom\'etries de la restriction de $(.,.)$ \`a un sous-espace non-d\'eg\'en\'er\'e $V_1$ de $V$;
\item $\mathfrak{g}_1$ est simple de rang r\'eel $\geq 2$; et
\item le groupe des isom\'etries de la restriction de $(.,.)$ \`a l'orthogonal de $V_1$ est compact.
\end{enumerate}
Il existe alors un sous-groupe de congruence $\Gamma \subset G$ tel que la repr\'esentation $\pi$ intervienne discr\`etement dans $L^2 (\Gamma \backslash G)$.
\end{thm}

\begin{conj} \label{relthc}
Soit $G$ comme au-dessus et soit $\pi = A_{\mathfrak{q}}$ une repr\'esentation cohomologique de $G$ telle que
$$\mathfrak{l}_0 \cong \mathfrak{l}_0' \oplus \mathfrak{g}_1$$
o\`u $ \mathfrak{l}_0' $ est une alg\`ebre de Lie compacte et $\mathfrak{g}_1$ est non-ab\'elienne du ``m\^eme type'' que $\mathfrak{g}_0$.
Alors, il existe un sous-groupe de congruence $\Gamma \subset G$ tel que la repr\'esentation $\pi$ intervienne discr\`etement dans $L^2 (\Gamma \backslash G)$.
\end{conj}

La Conjecture \ref{relthc} (qui est donc un corollaire de la Conjecture \ref{isolaut}) contient comme cas particuliers le Corollaire \ref{ttt} et le Th\'eor\`eme \ref{Li gen},
\`a titre de comparaison avec ce dernier, explicitons le Th\'eor\`eme \ref{relth} dans le cas des groupes unitaires, orthogonaux et symplectiques.

\begin{cor}
\begin{enumerate}
\item Soit $G=SU(p,q)$ avec $p\geq 2$. Alors la r\'eponse \`a la Question \ref{Q2} est positive pour chacune des repr\'esentations
$A((r^p),((q-s)^p))$, avec $r$, $s$ entiers naturels tels que $r+s \leq q-2$.
\item Soit $G=O(p,q)^0$ avec $p\geq 2$. Alors la r\'eponse \`a la Question \ref{Q1} est positive pour chacune des repr\'esentations
$A((r^p))$, avec $r$ entier naturel tel que $2r\leq \min (q-2 , p+q-5)$.
\item Soit $G=Sp(p,q)$. Alors la r\'eponse \`a la Question \ref{Q1} est positive pour chacune des repr\'esentations
$A((0), ((q-r)^p)_0$, avec $r$ entier naturel tel que $r\leq q$.
\end{enumerate}
\end{cor}

\medskip

Des restrictions simples sur la cohomologie des vari\'et\'es localement sym\'etriques d\'ecoulent de la classification des repr\'esentations cohomologiques~:
il est en effet imm\'ediat que si $X$ est irr\'eductible,
\begin{eqnarray} \label{annulation}
H^i (S(\Gamma )) = 0 , \  \mbox{ pour tout } 0<i< r_G ,
\end{eqnarray}
o\`u $r_G $ est l'infimum des ${\rm dim} (\mathfrak{u} \cap \mathfrak{p})$ sur toutes les sous-alg\`ebres paraboliques $\theta$-stable $\mathfrak{q} = \mathfrak{l}
+ \mathfrak{u}$ de $\mathfrak{g}$. On peut trouver une table des valeurs possibles de $r_G$ dans \cite[Table 8.2]{VoganZuckerman}. On a ainsi
$$r_G = \left\{
\begin{array}{ll}
\min (p,q) , & G= O(p,q)^0 ,\;  SU(p,q) \\
2 \min (p,q) , & G= Sp(p,q).
\end{array} \right.$$
En g\'en\'eral retenons simplement que $r_G$ est toujours inf\'erieur au rang r\'eel de $G$.

Rappelons que dans le cas hermitien l'existence d'une vari\'et\'e de Shimura $S(\Gamma )$ avec $H^{r_G } (S(\Gamma )) \neq 0$
est d\'emontr\'ee par Anderson \cite{Anderson}. Supposons $p\leq q$ et consid\'erons la restriction de $SU(p,q)$ \`a $O(p,q)^0$ (resp. de $SU(2p,2q)$ \`a $Sp (p,q)$), il d\'ecoule
de \cite{Lefschetz} (resp. de la m\'ethode de \cite{Lefschetz}) que l'unique repr\'esentation cohomologique du groupe $SU(p,q)$ (resp. $SU(2p,2q)$) de bidegr\'e $(p,0)$ (resp. $(2p,0)$)
contient discr\`etement l'unique repr\'esentation cohomologique de degr\'e $p$ (resp. $2p$) du groupe $O(p,q)^0$ (resp. $Sp(p,q)$), avec de plus une injection au niveau de la
$(\mathfrak{g} , K)$-cohomologie. Le point 1. de la Proposition \ref{P2} est donc v\'erifi\'e dans ces deux cas. Puisque par ailleurs la Proposition \ref{d} implique que l'unique
repr\'esentation cohomologique de degr\'e $p$ (resp. $2p$) du groupe $O(p,q)^0$ (resp. $Sp(p,q)$) est isol\'ee sous la condition $d=0$. Le Th\'eor\`eme d'Anderson (dans le cas $SU(p,q)$) et
la Proposition \ref{P2} permettent de compl\`eter l\'eg\`erement un r\'esultat de Li \cite{Li2} et d'en donner une d\'emonstration unifi\'ee.

\begin{cor}
Soit $G$ un groupe ${\Bbb Q}$-alg\'ebrique provenant (par restriction des scalaires) d'une alg\`ebre \`a involution
(de type I, II ou III) sur un corps de nombre et tel que $G^{\rm nc}=O(p,q)^0$ ou $Sp(p,q)$. Alors, il existe un sous-groupe de congruence
$\Gamma \subset G({\Bbb Q})$ tel que
$$H^{r_G} (S(\Gamma)) \neq 0 .$$
\end{cor}

L'analogue de ce r\'esultat est faux pour le groupe $U(p,q)$, cf. Proposition \ref{clozel} ou \cite{Clozel2}.

\section{Restriction de repr\'esentations unitaires}

Soit $G$ un groupe semi-simple r\'eel. Fixons $K$ un sous-groupe compact
maximal et $P_{\emptyset} = M_{\emptyset} A_{\emptyset} N_{\emptyset}$ un sous-groupe parabolique minimal de $G$. Rappelons qu'un sous-groupe parabolique $P=MAN
\subset G$ est dit {\it standard} s'il contient $P_{\emptyset}$ et qu'il est dit {\it
cuspidal} si ${\rm rang}_{{\Bbb C}} (M) = {\rm rang}_{{\Bbb C}} (K \cap M)$
(autrement dit si $M$ poss\`ede une s\'erie discr\`ete). Rappelons
maintenant la param\'etrisation de Langlands des repr\'esentations unitaires
de $G$, cf. \cite[Theorem 14.92]{Knapp}. (Nous adoptons ici pour les repr\'esentations induites et leurs quotients de Langlands des notations l\'eg\`erement diff\'erentes de celles du \S 1.)

\begin{thm} \label{LC}
Soit $P=MAN$ un sous-groupe parabolique cuspidal standard de $G$, soit
$\sigma$ une s\'erie discr\`ete de $M$ ou une limite non d\'eg\'en\'er\'ee
de s\'eries discr\`etes de $M$ et soit $\nu \in \mathfrak{a}^*$ avec ${\rm
Re} (\nu )$ dans la chambre de Weyl positive ouverte. La repr\'esentation
induite $\pi_{\sigma, \nu} = {\rm ind}_P^G (\sigma \otimes \nu )$ admet alors un unique
quotient irr\'eductible $J_{\sigma , \nu}$. Toute repr\'esentation
unitaire irr\'eductible de $G$ est soit temp\'er\'ee, soit de la forme
$J_{\sigma , \nu}$ avec $P$, $\sigma$ et $\nu$ comme ci-dessus.
\end{thm}

Dans cette section nous \'etudions la restriction des repr\'esentations unitaires
$J_{\sigma , \nu}$ \`a un sous-groupe $H$ de $G$ de {\bf rang (r\'eel)
\'egal \`a $\mathbf{1}$}. Consid\'erons donc $H \subset G$ un sous-groupe semi-simple
r\'eel de rang r\'eel \'egal \`a $1$ et tel que $K^H := H \cap K$ soit un sous-groupe compact maximal dans
$H$ (inclusion {\it standard}). Alors la restriction \`a $H$ de l'involution de Cartan $\theta$ de $G$
est une involution de Cartan et nous pouvons fixer un sous-groupe
parabolique minimal $P^H = M^H A^H N^H$ de $H$ tel que $A^H = A_{\emptyset} \cap H$ et
que la chambre de Weyl positive 
$$\mathfrak{a}_0^{H,+} \subset \mathfrak{a}_{\emptyset , 0}^{+}.$$
Le Th\'eor\`eme \ref{LC} s'applique \'evidemment \`a $H$, nous noterons $J^H_{\delta , \nu}$ ($\delta \in
\widehat{M^H}$) les repr\'esentations correspondantes. Notons enfin $\rho$
(resp. $\rho^H$) la demi-somme des racines positives de $\mathfrak{a}_{\emptyset}$
(resp. $\mathfrak{a}^H$) dans $\mathfrak{g}$ (resp. $\mathfrak{h}$).  

Le th\'eor\`eme
principal de cette section est alors~:

\begin{thm} \label{res}
Soient $P$, $\sigma$ et $\nu$ comme dans le Th\'eor\`eme \ref{LC} avec $J_{\sigma , \nu}$ {\bf unitaire}.
Supposons que la restriction de la repr\'esentation $J_{\sigma , \nu}$ au
groupe $H$ soit {\bf non temp\'er\'ee}.
Il existe alors un $M^H$-type $\delta$ contenu dans la restriction de
$\sigma$ \`a $M^H$ tel que  
la repr\'esentation $J^H_{\delta , (\nu -\rho)_{|\mathfrak{a}^H} + \rho^H}$
de $H$ soit faiblement contenue dans la restriction de $J_{\sigma , \nu}$
au groupe $H$. (Ici $\nu$ est \'etendue \`a tout $\mathfrak{a}_{\emptyset}$
de mani\`ere \`a \^etre \'egal \`a $0$ sur $\mathfrak{a}_M$.) 
\end{thm}

Avant de d\'emontrer le Th\'eor\`eme \ref{res}, commen\c{c}ons par quelques
rappels sur les fonctions radiales. Soit $\pi$ une repr\'esentation irr\'eductible unitaire de $G$ et $V$ le $(\mathfrak{g} ,K)$-module associ\'e.
Notons $I_{\pi }$ l'id\'eal du centre ${\cal Z} (\mathfrak{g})$ de l'alg\`ebre enveloppante d\'efini par
$$I_{\pi} = \{ Z \in  {\cal Z} (\mathfrak{g}) \; : \; \pi (Z) = 0 \mbox{ sur } V \}.$$
Rappelons (cf. \cite{CasselmanMilicic}) que l'id\'eal $I_{\pi}$ est cofini dans l'alg\`ebre ${\cal Z} (\mathfrak{g})$.

Soient $(\tau_i , V_{\tau_i})$, $i=1,2$, deux repr\'esentations irr\'eductibles de $K$ et $E={\cal L}(V_{\tau_2} , V_{\tau_1})$. 
Notons $\tau$ la repr\'esentation de $K \times K$ sur $E$~: $\tau (k_1 , k_2 ) \cdot L = \tau_1 (k_1 ) \circ L \circ \tau_2 (k_2^{-1} )$.
On appelle {\it fonction $\tau$-radiale} toute fonction
$$F : G \rightarrow E$$
v\'erifiant la condition de double $K$-\'equivariance suivante~:
\begin{eqnarray} \label{911}
F(k_1 g k_2 ) = \tau(k_1 , k_2^{-1} ) \cdot F(g) =  \tau_1 (k_1 ) F(g) \tau_2 (k_2 )
\end{eqnarray}
pour tous $g \in G$ et $k_1 , k_2 \in K$. Supposons maintenant que $\tau_1$ et $\tau_2$ sont deux $K$-types de la repr\'esentation $\pi$ et notons
$$i_{\tau_2} : V_{\tau_2} \hookrightarrow V$$
l'inclusion et
$$p_{\tau_1} : V \rightarrow V_{\tau_1}$$
la projection. Alors, la fonction $F: G \rightarrow E$ d\'efinie par 
\begin{eqnarray} \label{913}
F(g) = p_{\tau_1} \circ \pi (g) \circ i_{\tau_2} 
\end{eqnarray}
est $\tau$-radiale. Elle appartient en fait \`a ${\cal A} (G, \tau_1 , \tau_2 , I_{\pi} )$, l'ensemble des fonctions $\Phi$ lisses et $\tau$-radiales telles que
$R (Z) \Phi = 0$, pour tout $Z \in I_{\pi}$. (Avec $R$ repr\'esentation r\'eguli\`ere gauche dans l'espace des fonctions $\tau$-radiales.)
Une telle fonction est toujours analytique.

Int\'eressons-nous maintenant \`a l'asymptotique d'une fonction $\Phi \in
{\cal A} (G, \tau_1 , \tau_2 , I_{\pi} )$. D'apr\`es la d\'ecomposition de Cartan 
$$G= K \overline{A^+} K$$
il suffit de comprendre $\Phi (a)$ lorsque $a$ tend vers l'infini dans $\overline{A^+}$. Remarquons que la restriction de $\Phi$ \`a $A^+$ prend ses valeurs dans 
$$E^M := \{ L \in E \; : \; L= \tau (m,m) \cdot L = \tau_1 (m) L \tau_2 (m)^{-1} \mbox{ pour tout } m \in M \} , $$
o\`u $M$ d\'esigne le centralisateur de $A$ dans $K$.

Rappelons alors qu'il existe un ensemble fini $X \subset \mathfrak{a}^*$ et
un entier $d \in {\Bbb N}$ tels que toute fonction 
$\Phi \in {\cal A }(G, \tau_1 , \tau_2 , I_{\pi} )$ admet
un d\'eveloppement en s\'erie absolument convergente de la forme
\begin{eqnarray} \label{da}
\Phi (a) = \sum_{
\begin{array}{c}
\xi \in X - {\Bbb N} \Delta \\
|m| \leq d
\end{array} } (\log a)^m a^{\xi} c_{\xi , m} \;  \;  (a \in A^+ )
\end{eqnarray}
avec des coefficients uniquement d\'etermin\'es $c_{\xi , m } \in E^M$. (Ici $\Delta$ d\'esigne l'ensemble des racines simples de $\mathfrak{a}$ dans $\mathfrak{g}$.)
Un \'el\'ement $\xi \in X-{\Bbb N} \Delta$ tel qu'il existe un coefficient $c_{\xi , m} \neq 0$ pour un certain entier $m$ est appel\'e un {\it exposant} de $\Phi$ (le
long de $A^+$), nous notons ${\cal E}(\Phi )$ l'ensemble des exposants de $\Phi$. Le $K\times K$-type $\tau$ \'etant fix\'e, on appelle
{\it $\tau$-exposant} de $\pi$ tout \'el\'ement $\xi \in {\cal E} (\Phi )$
pour une certaine fonction lisse $\tau$-radiale $\Phi \in {\cal A}
(G,\tau_1 , \tau_2 , I_{\pi} )$.
Un tel exposant $\xi$ est dit {\it temp\'er\'e} si pour toute racine $\alpha \in \Delta$,
$$\langle {\rm Re} \xi + \rho , \omega_{\alpha} \rangle  \leq 0,$$
o\`u $\omega_{\alpha} \in \mathfrak{a}^*$ est le poids fondamental associ\'e \`a $\alpha$, {\it i.e.} $2 \langle \omega_{\alpha} , \beta \rangle / \langle \beta , \beta \rangle = \delta_{\alpha 
\beta}$.
.
La proposition suivante est bien connue, cf. \cite[Proposition 8.61]{Knapp}.
\begin{prop} \label{exp}
Soit $J_{\sigma, \nu}$ une repr\'esentation de $G$ comme dans le
Th\'eor\`eme \ref{LC}. Il existe alors un exposant de $J_{\sigma, \nu}$ de
la forme $\tilde{\nu} -\rho$ tel que $\tilde{\nu}_{|\mathfrak{a}} =
\nu$. 

Soient plus pr\'ecis\'ement $\tau_1$ et $\tau_2$ deux $K$-types
contenant (en restriction \`a $M\cap K$) un m\^eme $(M\cap K)$-type $\mu$
de la repr\'esentation $\sigma$. Il existe alors une fonction $\Phi \in
{\cal A} (G, \tau_1 , \tau_2 , I_{\pi})$ telle que 
\begin{enumerate}
\item il existe un exposant $\tilde{\nu} - \rho \in X$ (notations de
(\ref{da})) tel que $\tilde{\nu}_{|\mathfrak{a}} =
\nu$, et
\item il existe deux vecteurs $v_i \in V_{\tau_i}$, $i=1,2$, engendrant sous l'action
du groupe $M\cap K$ le $(M\cap K)$-type $\mu$ tels que $\langle
c_{\tilde{\nu} -\rho , 0} (v_2) , v_1 \rangle \neq 0$.
\end{enumerate}
\end{prop} 

\medskip

D\'etaillons le cas d'un groupe $H$ de rang $1$. Fixons
$\tau_1^H$ et $\tau_2^H$ deux $K^H$-types et consid\'erons $\delta$ une repr\'esentation
irr\'eductible de $M^H$ apparaissant dans chaque $(\tau_i^H )_{|K^H}$ ce que nous
notons~: 
$\delta \in \widehat{M^H} (\tau^H )$ (nous notons toujours $\tau^H =
(\tau_1^H , \tau_2^H)$). 
D\'esignons par $P_{\delta}^{\tau_1^H}$ le g\'en\'erateur de l'espace de
dimension $1$, ${\rm Hom}_{K^H} (\pi_{\delta,\nu}^H , V_{\tau_1^H })$ et 
par $J_{\delta}^{\tau_2^H }$ le g\'en\'erateur ($=(P_{\delta}^{\tau_2^H})^*$) de l'espace de dimension $1$, ${\rm Hom}_{K^H} ( V_{\tau_2^H} , \pi_{\delta ,\nu}^H )$. 
Si $\delta \in \widehat{M^H} (\tau^H)$ et $\nu \in (\mathfrak{a}^H)^*$, l'application 
$$\Phi_{\delta, \nu}^{\tau_1^H , \tau_2^H} : g \mapsto
P_{\delta}^{\tau_1^H} \circ \pi^H_{\delta, \nu} \circ  J_{\delta }^{\tau_2^H}$$
d\'efinit une fonction lisse, $\tau^H$-radiale. En fait (cf. Wallach \cite{Wallach}),
$${\cal A} (H, \tau_1^H , \tau_2^H , I_{\pi_{\delta , \nu}^H} ) = \langle
\Phi_{\delta , \nu}^{\tau_1^H , \tau_2^H } \; : \; \delta \in  \widehat{M^H} (\tau^H) \rangle .$$

Puisque $H$ est de rang $1$, il n'y a qu'une seule racine simple $\alpha$ de $\mathfrak{a}^H$ dans $\mathfrak{h}$. 
Notons $\omega = \omega_{\alpha}$. Dans la suite nous identifions le dual $(\mathfrak{a}^H)^*$ \`a ${\Bbb C}$ via l'isomorphisme $s \in {\Bbb C} \mapsto s \omega \in (\mathfrak{a}^H)^*$.
Rappelons alors le th\'eor\`eme classique suivant, cf. \cite{Wallach} ou \cite{Livre} pour une d\'emonstration.

\begin{thm} \label{W}
Fixons $\tau^H = (\tau_1^H , \tau_2^H) \in \widehat{K^H \times K^H}$. 
Alors, pour tout  $\delta \in \widehat{M^H} (\tau^H )$, il existe un r\'eel $\kappa>0$ tel que pour tout r\'eel strictement positif $\varepsilon$, il 
existe une constante $C_{\varepsilon}$ telle que pour tout \'el\'ement $\nu \in (\mathfrak{a}^H)^*$ v\'erifiant $\rho^H > {\rm Re} (\nu ) \geq \varepsilon$, il existe
des constantes $c_{\tau_1^H} (\nu )$ et $c_{\tau_2^H}  (\nu )$ telles que~:
$$|| \Phi_{\delta , \nu}^{\tau_1^H , \tau_2^H } (a) -  c_{\tau_2^H } (\nu )
a^{\nu - \rho^H } A_{\delta}
-  c_{\tau_1^H}  (\nu ) a^{- \rho^H - \overline{\nu}}  \tau_1 (k^* )
A_{\delta}^{\tau_1^H , \tau_2^H} \tau_2 (k^* )^{-1} || \leq C_{\varepsilon}
a^{-\rho^H -\kappa} ,$$  
o\`u $||.||$ d\'esigne la norme d'op\'erateur associ\'ee \`a des normes (quelconques) 
sur $V_{\tau_1^H}$ et $V_{\tau_2^H}$, $A_{\delta}^{\tau_1^H , \tau_2^H} = P_{\delta}^{\tau_1^H}
\circ J_{\delta}^{\tau_2^H} \in {\rm Hom}_{M^H} (V_{\tau_2^H} , V_{\tau_1^H})$,
et $k^*$ est un \'el\'ement de $K^H$ centralisant $A^H$ tel que $k^* e^{tX} (k^*)^{-1} = e^{-tX}$. 
\end{thm}

\medskip

\noindent
{\bf Remarque.} Les constantes $c_{\tau_1^H} (\nu )$ et $c_{\tau_2^H}  (\nu )$
sont explicites. Retenons simplement que l'on peut toutes deux les majorer
en modules et uniform\'ement par rapport au choix d'un $\nu$ v\'erifiant
$\rho^H  > {\rm Re} (\nu ) \geq \varepsilon$, par une constante
strictement positive $C(\varepsilon, \tau^H )$.

\medskip

Le Th\'eor\`eme \ref{LC} (et \cite[Theorem 16.10]{Knapp}) implique(nt) que toute repr\'esentation {\bf unitaire}
irr\'eductible {\bf non temp\'er\'ee} de $H$ est de la forme $J^H_{\delta , \nu}$ avec $\delta \in \widehat{M^H}$ et $\nu \in
(\mathfrak{a}^H)^*$ {\bf r\'eel} et tel que $\rho^H > \nu  >0$. 
Il d\'ecoule alors du Th\'eor\`eme \ref{W} qu'une telle repr\'esentation 
admet un unique exposant non temp\'er\'e $\nu - \rho^H
\in (\mathfrak{a}^H)^*$. 

Le Th\'eor\`eme \ref{W} permet de contr\^oler la croissance des coefficients des repr\'es- entations
irr\'eductibles de $H$. Commen\c{c}ons par rappeler \cite[Theorem
8.53]{Knapp} que si $\rho$ est une repr\'esentation {\bf temp\'er\'ee} de $H$
et $v_2$ deux vecteurs $K^H$-finis appartenant respectivement \`a deux
$K^H$-types $\tau_1^H$ et $\tau_2^H$, on a~:
\begin{eqnarray} \label{temp}
|\langle \rho (a) v_1 , v_2 \rangle | = O_{\tau^H} \left( ||v_1 || \cdot ||v_2 || \cdot
 a^{-\rho^H} \right),
\end{eqnarray}
o\`u les constantes implicites dans le $O_{\tau^H}$ sont uniformes \`a $\tau^H$ fix\'e.

Consid\'erons maintenant le cas d'une repr\'esentation
$J^H_{\delta , \nu}$ comme ci-dessus. Soient $v_1$
et $v_2$ deux vecteurs $K^H$-finis appartenant respectivement \`a deux
$K^H$-types $\tau_1^H$ et $\tau_2^H$. Fixons un r\'eel strictement positif $\varepsilon$ tel
que ${\rm Re} ( \nu ) > \varepsilon$. Alors~:
\begin{eqnarray} \label{nontemp}
\begin{array}{ccc}
\langle J^H_{\delta , \nu} (a) v_2 , v_1 \rangle & = & c_{\tau_2^H } (\nu )
\langle P_{\delta}^{\tau_2^H} v_2
, P_{\delta}^{\tau_1^H} v_1 \rangle_{V_{\delta}} \cdot a^{\nu-\rho^H} \\
& &  + O_{\varepsilon , \tau^H}
\left( ||v_1 || \cdot ||v_2 || \cdot a^{-\rho^H } \right),
\end{array}
\end{eqnarray}
o\`u chaque $P_{\delta}^{\tau_i^H}$ d\'esigne la projection $V_{\tau_i^H}
\rightarrow V_{\delta}$ et o\`u les constantes implicites dans le $O_{\varepsilon , \tau^H }$ sont
uniformes \`a $\varepsilon$ et $\tau^H$ fix\'es.

\medskip

L'\'etape principale de la d\'emonstration du Th\'eor\`eme \ref{res} est la
proposition suivante.

\begin{prop} \label{pgal}
Soit $\tau$ un $K$-type d'une repr\'esentation $J_{\sigma ,\nu}$ comme dans
le Th\'eor\`eme \ref{LC}. Soit $\xi\in \mathfrak{a}^*$ un $(\tau,
\tau)$-exposant de $J_{\sigma , \nu}$. Alors,
\begin{itemize}
\item soit la restriction $\xi_{| \mathfrak{a}^H } $ est temp\'er\'ee
comme \'el\'ement de $(\mathfrak{a}^H)^*$;
\item soit il existe un $M^H$-type $\delta$ apparaissant dans la
restriction de $\tau$ \`a $M^H$ tel que la repr\'esentation $J^H_{\delta ,
\nu '}$, avec $\nu ' = \xi_{| \mathfrak{a}^H } + \rho^H \in
(\mathfrak{a}^H)^*$, soit faiblement contenue dans la restriction de
$J_{\sigma ,\nu}$ \`a $H$.
\end{itemize}
\end{prop}
%en general on peut ne regarder que Re(xi_|) et dire que c'est un de Phi' a un multiple entier >0 pres (en ce sens il n'y a qu'un seul exposant dominant pour les induites 
{\it D\'emonstration.} La restriction de la repr\'esentation $\pi := J_{\sigma , \nu}$ au groupe $H$ se d\'ecompose en une somme directe~:
\begin{eqnarray} \label{dec en somme directe}
\pi_{|H} = \int_{\widehat{H}}^{\oplus} \rho d\mu (\rho ),
\end{eqnarray}
o\`u $\mu$ est une  mesure positive sur le dual unitaire $\widehat{H}$ de $H$. Son support est par d\'efinition le support de la restriction
de $\pi$ \`a $H$.
Fixons $\xi \in \mathfrak{a}^*$ un $(\tau, \tau)$-exposant de $\pi$ comme dans l'\'enonc\'e de la Proposition \ref{pgal}
et notons $\xi '$ l'\'el\'ement $\xi_{|\mathfrak{a}^{H} } \in
(\mathfrak{a}^H )^*$ que nous supposerons non temp\'er\'e. Remarquons alors imm\'ediatement que la d\'emonstration de \cite[Theorem 16.10]{Knapp} implique
que $\xi '$ est {\bf r\'eel}.

Nous d\'eduirons facilement la Proposition \ref{pgal} du lemme suivant.

\begin{lem} \label{lgal}
Il existe une repr\'esentation $\rho$ de $H$ appartenant
au support de $\mu$ et poss\`edant deux $K^H$-types $\tau_1^H$ et $\tau_2^H \subset \tau_{|K^H}$ telle que
$\xi '$ soit un $\tau^H$-exposant de $\rho$ (o\`u $\tau^H = (\tau_1^H ,
\tau_2^H )$).
\end{lem}
{\it D\'emonstration.} Afin de ne pas allourdir le texte nous notons dans la suite $(P (\xi ') )$ la
conclusion du Lemme \ref{lgal} relativement \`a un \'el\'ement $\xi' \in (\mathfrak{a}^H )^*$.

Fixons deux vecteurs unitaires $K$-finis $v_1$ et $v_2$ appartenant
\`a l'espace de la repr\'esentation $\tau$ et tels qu'il existe une
fonction $\Phi \in {\cal A} (G, \tau , \tau , I_{\pi})$ et un entier $m$
tels que l'exposant $\xi \in {\cal E} (\Phi)$ (notations de (\ref{da})) et le coefficient 
$\langle c_{\xi , m} v_2 ,v_1 \rangle$ soit non nul.
Remarquons que l'on peut toujours supposer $v_1$ (resp. $v_2$) appartenant \`a l'espace 
d'une repr\'esentation $\tau_1^H$ (resp. $\tau_2^H$) de $K^H$ (apparaissant dans la
restriction du $K$-type $\tau$).

La d\'ecomposition (\ref{dec en somme directe}) implique la d\'ecomposition suivante~:
\begin{eqnarray} \label{dec des coeff}
\langle \pi (a) v_2 , v_1 \rangle = \int_{\widehat{H} } \langle \rho (a)
v_{2, \rho} , v_{1, \rho}  \rangle d\mu ( \rho ) ,
\end{eqnarray}
o\`u $v_i  = \int_{\widehat{H}}^{\oplus} v_{i ,\rho}$, $i=1,2$ et $a$ est un \'el\'ement
de $A^H$. Remarquons que dans la d\'ecomposition (\ref{dec des coeff}) on
peut supposer que la mesure $\mu$ ne charge que les repr\'esentations
contenant les deux $K^H$-types $\tau_1^H$ et $\tau_2^H$.

Le groupe $H$ est de rang $1$. Chaque repr\'esentation $\rho$ dans
(\ref{dec des coeff}) a au plus un exposant non temp\'er\'e
$\xi_{\rho} \in (\mathfrak{a}^H)^*$. 

\begin{lem} \label{llll}
Soit $\xi ' \in (\mathfrak{a}^H)^*$ non temp\'er\'e tel que la propri\'et\'e $(P(\xi '))$ soit viol\'ee. Il existe alors un voisinage ouvert $\Omega$  
de $\xi '$ dans $(\mathfrak{a}^H)^*$ tel que la propri\'et\'e $(P(\alpha ))$ soit viol\'ee pour tout $\alpha \in \Omega$.
\end{lem}
{\it D\'emonstration.} Nous montrons la contrapos\'ee. Supposons donc que pour tout voisinage ouvert $\Omega$ de 
$\xi '$ dans $(\mathfrak{a}^H)^*$, il existe un \'el\'ement $\alpha \in \Omega$ tel que la propri\'et\'e $(P(\alpha ))$ soit v\'erifi\'ee.
Il existe alors une repr\'esentation $\rho_{\alpha}$ de $H$ appartenant au
support de $\mu$ et poss\`edant deux $K^H$-types $\tau_{\alpha,1}^H$ et
$\tau_{\alpha,2}^H \subset \tau_{|K^H}$ telle que $\alpha$ soit un
$\tau_{\alpha}^H$-exposant  de $\rho_{\alpha}$.

Fixons $\Omega_j$ une suite de voisinages ouverts de $\xi '$ d'intersection
r\'eduite \`a $\xi '$ et $\alpha_j \in \Omega_j$ un \'el\'ement non temp\'er\'e tel que $(P(\alpha_j ))$ soit v\'erifi\'ee.
D'apr\`es le Th\'eor\`eme \ref{LC}, chaque repr\'esentation (non
temp\'er\'ee) $\rho_j := \rho_{\alpha_j}$ est \'equivalente \`a un quotient
de Langlands $J^H_{\delta_j , \nu_j}$ avec $\delta_j \in \widehat{M^H}$,
$\nu_j \in \mathfrak{a}_H^*$, ${\rm Re} (\nu_j )>0$ et $\alpha_j = \nu_j -\rho^H$.

Nous allons d\'emontrer qu'il existe une sous-suite de $(\rho_j )$ qui
converge dans $\widehat{H}$ vers une repr\'esentation $\rho$ poss\`dant
deux $K^H$-types $\tau_1^H$ et $\tau_2^H \subset \tau_{|K^H}$ telle que
$\xi'$ soit un $\tau^H$-exposant de $\rho$.

La restriction de $\tau$ au compact $K^H$ ne contient qu'un nombre fini de sous-repr\'esentations irr\'eductibles.
Quitte \`a extraire une sous-suite de $(\rho_j )$, nous pouvons donc
supposer que tous les $K^H \times K^H$-types $\tau_{\alpha_j}^H$ sont
\'egaux \`a un certain $K^H \times K^H$-type $\tau^H = (\tau_1^H , \tau_2^H)$.
Remarquons maintenant que puisque $H$ est de rang $1$, le groupe $M^H$ est
compact. Puisque la repr\'esentation $\rho_j$ poss\`ede $\tau_1^H$ (resp. $\tau_2^H$) comme $K^H$-type,
le th\'eor\`eme de r\'eciprocit\'e de Frobenius implique que la
repr\'esentation $\delta_j$ est contenue dans la restriction de
$\tau_1^H$ (resp. $\tau_2^H$) \`a $M^H$. Ces repr\'esentations
ne contiennent, l\`a encore, qu'un nombre fini de sous-repr\'esentations
irr\'eductibles. Quitte \`a extraire encore une fois, nous pouvons donc supposer tous les $\delta_i$ \'egaux
\`a une m\^eme repr\'esentation $\delta \in \widehat{M^H}$.

La suite des exposants $\alpha_j = \nu_j -\rho^H \in (\mathfrak{a}^H)^*$
converge vers $\xi '$ dans $(\mathfrak{a}^H)^*$, la suite $(\nu_j )$ converge donc vers $\nu = \rho^H + \xi '$. Et la suite de repr\'esentations $(\rho_j )$
converge dans $\widehat{H}$ vers une repr\'esentation $\rho$ (obtenue comme sous-quotient de l'induite
${\rm ind}_{P^H}^H (\delta \otimes \nu )$) contenant les $K^H$-types 
$\tau_1^H$ et $\tau_2^H$ et telle que $\xi'$ soit un $\tau^H$-exposant de $\rho$.

Le support de $\mu$ \'etant ferm\'e dans $\widehat{H}$, la repr\'esentation
$\rho \in \mbox{support } \mu$ et la propri\'et\'e $(P(\xi' ))$ est v\'erifi\'ee. Ce qui conclut la d\'emonstration
(par contrapos\'ee) du Lemme \ref{llll}.

\bigskip

Chacun des coefficients $\langle \rho (a) v_{2, \rho} , v_{1, \rho}  \rangle$ 
dans (\ref{dec des coeff}) est une fonction $\tau^H$-radiale du groupe $H$
de rang $1$. On peut donc leur appliquer les estim\'ees (\ref{temp}) et
(\ref{nontemp}) pour un certain r\'eel strictement positif fix\'e $\varepsilon <
\xi' +\rho^H $ ($\xi'$ est non temp\'er\'e). 
On a alors~:
\begin{eqnarray} \label{dec2}
\begin{array}{lcl} 
\langle \pi (a) v_2 , v_1 \rangle&  =  & \int_{\rho \in \widehat{H} , \; \xi_{\rho} + \rho^H > \varepsilon} c_{\tau_2}
(\xi_{\rho}+\rho^H)  \langle P_{\delta_{\rho}} v_{2 , \rho} ,
P_{\delta_{\rho}} v_{1 , \rho} \rangle a^{\xi_{\rho}}  d\mu ( \rho ) \\
& &  +
O_{\varepsilon , v_1  , v_2 } (a^{-\rho^H} ) .
\end{array}
\end{eqnarray}
Chaque $\xi_{\rho}$ est par ailleurs r\'eel (chaque repr\'esentation $\rho$ est unitaire non temp\'er\'ee et $H$ est de rang $1$) et (Lemme \ref{llll}) il existe un r\'eel $\eta>0$ tel que   
$\xi_{\rho} \notin ]\xi ' -\eta , \xi' + \eta[$. Il d\'ecoule alors de (\ref{dec2}) que
$$\langle \pi (a) v_2 , v_1 \rangle = \frac{1}{O_{\varepsilon,  v_1 , v_2} ( a^{-\xi ' - \eta \omega})} + O_{\varepsilon , v_1 , v_2} (a^{\xi ' -\eta \omega}) ,$$
pour tout $a \in A^{H,+}$. Ce qui contredit le fait que $\xi$ soit un exposant de $\langle \pi (a) v_2 , v_1 \rangle$ et
conclut le d\'emonstration du Lemme \ref{lgal}.

\bigskip

D\'emontrons maintenant la Proposition \ref{pgal}. Supposons la restriction $\xi_{| \mathfrak{a}^H } $ non temp\'er\'ee
comme \'el\'ement de $(\mathfrak{a}^H)^*$. D'apr\`es le Lemme \ref{pgal}, la restriction de la repr\'esentation $\pi$ au sous-groupe $H$
contient alors faiblement une repr\'esentation $\rho$ poss\`edant deux $K^H$-types $\tau_1^H$ et $\tau_2^H \subset \tau_{|K^H}$ et telle que
$\xi '$ soit un $\tau^H$-exposant de $\rho$. La repr\'esentation est alors n\'ecessairement de la forme $J_{\delta , \nu}^H$ avec 
$\delta \in \widehat{M^H} (\tau^H )$ et $\nu ' = \xi' + \rho^H$. Et la Proposition \ref{pgal} est d\'emontr\'ee.

\bigskip

Le Th\'eor\`eme \ref{res} est une cons\'equence imm\'ediate de la Proposition \ref{exp} et de la d\'emonstration de la Proposition \ref{pgal}. 

\bigskip

\noindent
{\bf Remarques.} {\bf 1.} Testons la Proposition \ref{pgal}, dans 
le cas simple du groupe $H=SL_2 ({\Bbb R})$ plong\'e diagonalement dans $G= SL_2 ({\Bbb R}) \times SL_2 ({\Bbb R})$. Les repr\'esentations non temp\'er\'ees appartiennent
\`a la s\'erie compl\'ementaire. Soient donc $\pi_s$ et $\pi_r$ ($r, s \in ]0,1[$) deux repr\'esentations de la s\'erie compl\'ementaire. Le produit tensoriel 
$\pi_s \otimes \pi_r$ est une repr\'esentation de $G$, elle admet une fonction radiale dont un coefficient dominant est \'egal \`a 
$$\frac{s-1}{2} \varepsilon_1 + \frac{r-1}{2} \varepsilon_2 .$$
En restriction au sous-groupe diagonal de $H$, celui-ci est \'egal \`a $(r+s-2)/2$, il est temp\'er\'e si et seulement si $r+s\leq1$. Dans le cas contraire 
$r+s > 1$, l'\'el\'ement $(r+s-2)/2$ doit intervenir comme coefficient d'une repr\'esentation faiblement contenue dans la restriction \`a $H$ du produit
tensoriel $\pi_s \otimes \pi_r$. C'est effectivement le cas puisqu'alors un th\'eor\`eme classique de Repka \cite{Repka} affirme que la repr\'esentation 
$\pi_{r+s-1}$ intervient discr\`etement dans la restriction \`a $H$ de la repr\'esentation $\pi_s \otimes \pi_r$ de $G$. %pe ecrire le thm en termes de restriction de quotient de L et de
%faiblt contenue ca serait plus parlant.: dire si pi a un K type + un coeff ... alors la res contient faiblement l'unique sous-quotient de L tq ...
%cas interessant : celui des rep cohom puisqu'avec isolation ca permet de montrer thm

{\bf 2.} La Proposition \ref{pgal} implique imm\'ediatement le principe \'enonc\'e dans la Remarque 4 apr\`es le Corollaire 1.3 de \cite{BurgerSarnak}.

\bigskip

Soit $G$ un groupe semi-simple r\'eel et soit $H \subset G$ un sous-groupe semi-simple r\'eel (pas n\'ecessairement de rang $1$) tel que
$K^H :=H\cap K$ soit un sous-groupe compact maximal dans $H$ et toujours $(\mathfrak{a}^H_0 )^+ \subset \mathfrak{a}_0^+$. Nous conjecturons
alors la g\'en\'eralisation suivante de la Proposition \ref{pgal}.

\begin{conj} \label{cgal}
Soient $\pi \in \widehat{G}$, $\tau$ un $K$-type de $\pi$ et $\xi$ un $\tau$-exposant de $\pi$. Il existe alors
une repr\'esentation $\sigma \in \widehat{H}$ faiblement contenue dans la restriction de $\pi$ \`a $H$, un
$K^H$-type $\tau^H \subset \tau_{|K^H}$ de $\sigma$ et un $\tau^H$-exposant $\nu$ de $\sigma$ tels que
$${\rm Re} (\nu ) = \max \left( {\rm Re} (\xi)_{|\mathfrak{a}^H} - \rho_{|\mathfrak{a}^H} + \rho^H , 0 \right).$$
\end{conj}

Cette Conjecture est peut-\^etre accessible via la th\'eorie de Mackey telle que d\'evelopp\'ee par Venkatesh dans \cite{Venkatesh}, lorsque $H$ est un {\bf sous-groupe de Levi} de $G$ et $\xi$ un exposant dominant. 
Tel que nous comprenons \cite{Venkatesh} le cas du groupe $G=GL(n)$ (et $H= \prod GL(a_i )$ Levi de $G$),
dont on connait explicitement le dual unitaire, semble en d\'ecouler.

Contentons-nous ici de remarquer que la Conjecture \ref{cgal} est naturellement reli\'ee \`a la combinatoire d\'ecrite
par Clozel dans \cite{Clozel5} (cf. aussi \cite{Venkatesh}). Comme dans \cite{Clozel5} et \cite{Venkatesh}, \'etudions
le cas du groupe $GL(n)$. Celui-ci est particuli\`erement agr\'eable car nous connaissons son dual unitaire \cite{Vogan4}~: soient
$m$ un entier $=1$ ou $2$ et $\delta$ une repr\'esentation de la s\'erie discr\`ete de $GL(m,{\Bbb R})$. Pour tout entier naturel $j$, la repr\'esentation de
$GL(mj,{\Bbb R})$ induite de $(\delta |{\rm det}|^{(j-1)/2} \otimes \delta |{\rm det}|^{(j-3)/2} \otimes \ldots \otimes \delta |{\rm det}|^{(1-j)/2})$ admet un unique
quotient irr\'eductible~: $u(\delta , j)$. Si $0<\alpha <1/2$, la repr\'esentation de $GL(2mj,{\Bbb R})$ induite de $u(\delta , j)|{\rm det}|^{\alpha} \otimes
u(\delta , j) |{\rm det}|^{-\alpha}$ est unitarisable; nous la notons $u(\delta , j)[\alpha , -\alpha]$.

Toute repr\'esentation unitaire de $GL(n,{\Bbb R})$ est induite (unitaire) de repr\'esenta- tions du type $u(\delta , j)$ et $u(\delta , j)[\alpha , -\alpha]$,
et cette expression est unique \`a permutation pr\`es.
\`A une telle repr\'esentation unitaire $\pi$ nous associons une matrice diagonale constitu\'ee des blocs
$$\left(
\begin{array}{cccc}
(j-1)/2 & & & \\
& (j-3)/2 & &  \\
& & \ddots & \\
& & & (1-j) /2
\end{array} \right) $$
pour chaque $u(\delta, j)$ et r\'ep\'et\'e $m$ fois, et des deux blocs
$$\left(
\begin{array}{cccc}
(j-1)/2 \pm \alpha & & & \\
& (j-3)/2 \pm \alpha & &  \\
& & \ddots & \\
& & & (1-j) /2 \pm \alpha
\end{array} \right) $$
pour chaque $u(\delta , j) [\alpha , -\alpha]$ et r\'ep\'et\'es $m$
fois. Nous notons $T_{\pi}$ cette matrice r\'eordonn\'ee de telle mani\`ere
que
$$T_{\pi} = \left(
\begin{array}{ccc}
T_1 & & \\
 & \ddots & \\
& & T_n
\end{array} \right) , \; T_1 \geq \ldots  \geq T_n.$$
Soit $A \subset GL(n,{\Bbb R})$ le tore diagonal et soit
$$A_+ = \left\{ \left(
\begin{array}{ccc}
x_1 & & \\
 & \ddots & \\
& & x_n
\end{array} \right) \in H \; : \; x_1 \geq \ldots \geq x_n \geq 0
 \right\}.$$

Soit $H=GL(m,{\Bbb R})$ plong\'e de mani\`ere standard dans $G=GL(n,{\Bbb
R})$ (de telle mani\`ere que $A^{H,+} \subset A^+$).
La Conjecture \ref{cgal} implique alors que la restriction de $\pi$ au groupe
$H$ contient faiblement une repr\'esentation $\sigma$ telle que
$$T_{\sigma} = \langle T_{\pi} - {\rm diag}(\frac{n-1}{2} , \ldots , \frac{1-n}{2} ) + {\rm diag} (
\frac{m-1}{2} , \ldots , 0 , \ldots , 0 , \ldots ,\frac{1-m}{2} ) \rangle_m ,$$
o\`u la notation $\langle \cdot , \cdot \rangle_m$ signifie que nous ne gardons que les
termes appartenant \`a $GL(m)$ en rempla\c{c}ant ceux d'entre eux qui sont n\'egatifs par $0$.

Nous retrouvons en particulier la combinatoire sur les $SL(2)$-types d\'ecrite dans \cite{Venkatesh}.

\bigskip

Revenons maintenant \`a la Conjecture \ref{isolaut}. Nous avons rappel\'e que la Conjecture \ref{ArthurC} implique qu'un
groupe $L_0$ (localement isomorphe \`a) $SO(n,1)$ ($n\geq 2$) ou $SU(n,1)$ ($n\geq 1$) doit v\'erifier l'hypoth\`ese
\begin{eqnarray*}
({\rm Hyp}) & : & \left\{
\begin{array}{l}
\mbox{Si } \pi = J_{\delta , s} \in (\widehat{L_0} )_{{\rm Aut}} \mbox{ est non temp\'er\'ee, alors} \\
s \in \frac12 {\Bbb Z} \; (0 < s \leq \rho^{L_0} ) .
\end{array}
\right.
\end{eqnarray*}
Cette hypoth\`ese est naturellement approch\'ee par les hypoth\`eses~:
\begin{eqnarray*}
({\rm Hyp}_{\varepsilon}) & : & \left\{
\begin{array}{l}
\mbox{Si } \pi = J_{\delta , s} \in (\widehat{L_0} )_{{\rm Aut}} \mbox{ est non temp\'er\'ee, alors} \\
\mbox{ soit } s \in \frac12 {\Bbb Z} \; (0 < s \leq \rho^{L_0}); \;  \mbox{ soit } 0<s\leq \varepsilon ,
\end{array}
\right.
\end{eqnarray*}
avec $\varepsilon >0$. Toute hypoth\`ese $({\rm Hyp}_{\varepsilon} )$ avec
$\varepsilon < \rho^{L_0}$ est non triviale.

Dans \cite{Livre}, avec Clozel, nous montrons l'hypoth\`ese (${\rm Hyp}_{\frac{4}{5}}$) pour le groupe
$SU(2,1)$.

Fixons maintenant $\mathfrak{q}$ une sous-alg\`ebre parabolique $\theta$-stable de $\mathfrak{g}$ et $L \subset G$ le
sous-groupe de Levi associ\'e. Selon la Conjecture \ref{conj gal}, la repr\'esentation $A_{\mathfrak{q}}$ ne devrait
appara\^{i}tre dans le dual automorphe $\widehat{G}_{{\rm Aut}}$ que si (quitte \`a conjuguer $L$ dans $G$) $L$ est
un sous-groupe rationnel de $G$. Nous le supposerons par la suite.
Supposons enfin que $L$ poss\`ede un facteur simple non compact $L_0$, localement isomorphe \`a $SO(n,1)$ ($n\geq 2$)
ou $SU(n,1)$ ($n\geq 1$). La repr\'esentation $\pi :=A_{\mathfrak{q}}$ est alors non-isol\'ee dans le dual
unitaire $\widehat{G}$ de $G$. Le Th\'eor\`eme \ref{res} (ou plut\^ot la Proposition \ref{pgal}) permet n\'eanmoins de montrer~:

\begin{thm} \label{rel}
Supposons que le groupe $L_0$ v\'erifie l'hypoth\`ese $({\rm Hyp}_{\varepsilon})$. Alors, la
repr\'esentation $\pi$ est {\bf isol\'ee} dans
$$\{ \pi \} \cup \widehat{G}_{{\rm Aut}}$$
d\`es que $2\rho^{L_0} - \rho_{|\mathfrak{a}^{L_0}} > \varepsilon$.
\end{thm}
{\it D\'emonstration.} Nous montrons la contrapos\'ee. Supposons donc
$2\rho^{L_0} - \rho_{|\mathfrak{a}^{L_0}} > \varepsilon$ et qu'il existe une suite $(\pi_i )$ de repr\'esentations dans $\widehat{G}_{{\rm Aut}}$
qui converge vers $\pi$, avec $\pi_i \not\cong \pi$.
Il d\'ecoule de la description des repr\'esentations cohomologiques faite au \S 1 que 
$\pi$ s'obtient comme quotient de Langlands d'une induite unitaire \`a partir d'un sous-groupe parabolique cuspidal $P=MAN$, 
o\`u $A$ est un sous-groupe de Cartan maximalement d\'eploy\'e dans $L$, et $\pi = J_{\sigma , \rho^L}$. La Proposition \ref{exp} implique alors
que la repr\'esentation $\pi$ admet un exposant dont la restriction \`a $\mathfrak{a}^{L_0}$ est \'egale \`a $\rho^{L_0} - \rho_{|\mathfrak{a}^{L_0}} \in
\frac12 {\Bbb Z}$ et $>\varepsilon$. Pour $i$
suffisamment grand, la repr\'esentation $\pi_i$ admet donc un exposant dont la restriction \`a $\mathfrak{a}^{L_0}$ est $>\varepsilon - \rho^{L_0}$ et $\notin \frac12 {\Bbb Z}$.
Et la Proposition \ref{pgal} implique que le support de la restriction de $\pi_i$ au groupe $L_0$ contient une repr\'esentation de la forme $J_{\delta ,s}$ avec $s >\varepsilon$ et $s \notin \frac12 {\Bbb Z}$.
Mais puisque le sous-groupe $L_0$ est rationnel, la fonctorialit\'e de la restriction (Burger et Sarnak, cf. \S 4) implique que le support de la restriction de $\pi_i$ au groupe $L_0$
est contenu dans $(\widehat{L_0})_{{\rm Aut}}$, ce qui contredit l'hypoth\`ese $({\rm Hyp}_{\varepsilon})$ et conclut la d\'emonstration du Th\'eor\`eme \ref{rel}.

\bigskip

Admettons temporairement qu'un groupe du type $SU(n,1)$ v\'erifie toujours l'hypoth\`ese (Hyp).
Le Th\'eor\`eme \ref{rel} s'applique alors aux repr\'esentations $A((a^p) , ((a+1)^p))$, $0 \leq a \leq q-1$,
du groupe $SU(p,q)$ lorsque $p\geq q$. N\'eanmoins le Th\'eor\`eme \ref{res} n'est pas employ\'e ici de mani\`ere optimale~:
il n'y a {\it a priori} aucune raison de restreindre $A_{\mathfrak{q}}$ au groupe $L_0$, nous pourrions tout aussi
bien restreindre la repr\'esentation cohomologique \`a un sous-groupe rationnel $H$ de rang $1$ diff\'erent de $L_0$. En choisissant
$H=SU(1,q)$, on peut ainsi v\'erifier qu'une repr\'esentation cohomologique
$\pi = A((a^p) , ((a+1)^p))$, avec $0 \leq a \leq q-1$, du groupe $G=SU(p,q)$ ``standard'' (c'est-\`a-dire provenant d'un vrai groupe unitaire sur un corps de nombre)
devrait (toujours modulo l'hypoth\`ese (Hyp) pour $SU(n,1)$) toujours \^etre {\bf isol\'ee} dans
$$\{ \pi \} \cup \widehat{G}_{{\rm Aut}}.$$
(On a en effet dans ce cas $\rho^L = p$, $\rho^H = q$ et $\rho_{|\mathfrak{a}^H} = p+q-1$.)
De la m\^eme mani\`ere on peut \'egalement d\'emontrer la l\'eg\`ere g\'en\'eralisation suivante du Th\'eor\`eme 2 de \cite{Livre}~:

\begin{thm} \label{livre}
Soit $H\subset G$ une inclusion standard entre deux ${\Bbb Q}$-groupes semi-simples de rang r\'eel $1$. Supposons que le groupe $H$ v\'erifie
l'hypoth\`ese $({\rm Hyp}_{\varepsilon})$, pour un certain r\'eel positif $\varepsilon$, alors le groupe $G$ v\'erifie l'hypoth\`ese
$({\rm Hyp}_{\rho^G -\rho^H + \varepsilon} )$.
\end{thm}

Puisque le groupe $SU(2,1)$ v\'erifie l'hypoth\`ese (${\rm Hyp}_{\frac{4}{5}}$), il d\'ecoule du Th\'eor\`eme \ref{livre} qu'un groupe $SU(n,1)$ ``standard'' v\'erifie l'hypoth\`ese $({\rm Hyp}_{n-6/5})$.
Les rep- r\'esentations cohomologiques de degr\'e fortement primitif $1$ sont donc isol\'ees dans le dual automorphe.
Dans le cas du groupe $SU(n,2)$, il d\'ecoule imm\'ediatement
du Corollaire \ref{isolU} qu'une repr\'esentation cohomologique de degr\'e fortement primitif $\leq n-1$ est toujours isol\'ee dans le dual unitaire. La
repr\'esentation cohomologique de $SU(n,2)$ de bidegr\'e fortement primitif $(n,0)$~:
$$A((1^n) , (2^n))$$
est quant \`a elle non isol\'ee dans le dual unitaire du groupe $SU(n,2)$.
Puisque $2n -(n+2-1) = n-1> n- 6/5$, l'hypoth\`ese (${\rm Hyp}_{n-6/5}$) pour le groupe $SU(n,1)$ ``standard'' et le Th\'eor\`eme \ref{rel} impliquent inconditionnellement~:

\begin{thm} \label{i}
Soit $G$ un groupe alg\'ebrique sur ${\Bbb Q}$ obtenu, par restriction des scalaires, \`a partir d'un groupe unitaire sur un corps de nombres totalement r\'eel ({\it vrai groupe unitaire}) et tel
que $G^{{\rm nc}} = SU(n,2)$. Alors la repr\'esentation cohomologique holomorphe $\pi$ de $G$ de degr\'e fortement primitif $n$ est {\bf isol\'ee} dans 
$$\{ \pi \} \cup \widehat{G}_{{\rm Aut}}.$$
\end{thm}

\bibliography{bibliographie}

\bibliographystyle{plain}

\bigskip

\noindent
Unit\'e Mixte de Recherche 8553 du CNRS, \\
D\'epartement de math\'ematiques et applications (DMA), \\ 
45, rue d'Ulm 75230 Paris Cedex 05, France \\
{\it adresse electronique :} \texttt{Nicolas.Bergeron@ens.fr}

\end{document}